\documentclass[12pt,reqno]{amsart}
\usepackage{amsmath}
\usepackage{amssymb}
\usepackage{amsfonts}
\numberwithin{equation}{section}
\newcommand{\D}{\mathbb{D}\,}
\newcommand{\T}{\mathbb{T}\,}
\newcommand{\Z}{\mathbb{Z}\,}
\newcommand{\NN}{\mathbb{N}\,}
\newcommand{\CC}{\mathbb{C}\,}

\newcommand{\ahd}{{\mathcal{A}_h(\D)}}
\newcommand{\ahc}{{\mathcal{A}_h(\CC)}}
\newcommand{\ahdp}{{\mathcal{A}^p_h(\D)}}
\newcommand{\ahdd}{{\mathcal{A}^2_h(\D)}}
\newcommand{\ahcd}{{\mathcal{A}^2_h(\CC)}}
\newcommand{\ahcp}{{\mathcal{A}^p_h(\CC)}}
\newcommand{\RE}{\mathop{\rm Re\,}}
\newcommand{\IM}{\mathop{\rm Im\,}}
\newcommand{\hol}{\mathop{\rm Hol\,}}
\newcommand{\ve}{\varepsilon}
\newcommand{\dmGD}{D^{-}_\rho(\Gamma,\D)}
\newcommand{\dpGD}{D^{+}_\rho(\Gamma,\D)}
\newcommand{\dmGC}{D^{-}_\rho(\Gamma,\CC)}
\newcommand{\dpGC}{D^{+}_\rho(\Gamma,\CC)}
\newcommand{\conod}{{($\rm I_\D$\!)}}
\newcommand{\condd}{{($\rm II_\D$\!)}}
\newcommand{\conoc}{{($\rm I_\CC$\!)}}
\newcommand{\condc}{{($\rm II_\CC$\!)}}

\DeclareMathOperator{\Aut}{Aut}
\DeclareMathOperator{\dist}{dist}

\DeclareMathOperator{\card}{Card}
\DeclareMathOperator{\clos}{clos}
\DeclareMathSymbol{\subsetneqq}{\mathbin}{AMSb}{36}
\DeclareMathOperator{\ind}{ind}

\newtheorem{th1}{{\bf Theorem}}[section]

\newtheorem{thm}[th1]{{\bf Theorem}}
\newtheorem{lem}[th1]{{\bf Lemma}}

\newtheorem{prop}[th1]{{\bf Proposition}}
\newtheorem{cor}[th1]{{\bf Corollary}}

\begin{document} 
\author[A. Borichev, R. Dhuez, K. Kellay]
{A. Borichev, R. Dhuez, K. Kellay}

\title[Sampling and interpolation in radial weighted spaces]
{Sampling and interpolation in radial weighted spaces of analytic
functions}

\begin{abstract} We obtain sampling and interpolation theorems 
in radial weighted spaces of analytic functions for weights of arbitrary 
(more rapid than polynomial) growth. We give an application to invariant subspaces
of arbitrary index in large weighted Bergman spaces.
\end{abstract}

\date{\today}

\maketitle

\section{Introduction}

Let $h:[0,1)\to[0,+\infty)$ be an increasing function such that
$h(0)=0$, and $\lim_{r\to1}h(r)=+\infty$. We extend $h$ by
$h(z)=h(|z|)$, $z\in \D$, and call such $h$ a weight function.
Denote by $\ahd$ the Banach space of holomorphic functions on the
unit disk $\D$ with the norm
$$
\|f\|_h=\sup_{z\in\D}|f(z)|e^{-h(z)}<+\infty.
$$

A subset $\Gamma$ of $\D$ is called a {\it sampling set} for $\ahd$
if there exists $\delta>0$ such that for every $f\in\ahd$ we have
$$
\delta \|f\|_h\le
\|f\|_{h,\Gamma}=\sup_{z\in\Gamma}|f(z)|e^{-h(z)}.
$$

A subset $\Gamma$ of $\D$ is called an {\it interpolation set}
for $\ahd$ if for every function $a$ defined on $\Gamma$ such that $\|a\|_{h,\Gamma}<\infty$
there exists $f\in \ahd$ such that
$$
a=f\bigm |  \Gamma.
$$
In this case there exists $\delta=\delta(h,\Gamma)>0$ such that for every $a$
with $\|a\|_{h,\Gamma}<\infty$ we can find such $f\in \ahd$ with
$$
\delta \|f\|_h\le \|a\|_{h,\Gamma}.
$$

Next we assume that $h\in C^2(\D)$, and 
$$
\Delta h(z)=
\Bigl(\frac{\partial^2}{\partial x^2}+\frac{\partial^2}{\partial y^2}\Bigr)h(x+iy)\ge 1, 
\qquad z\in \D.
$$

We consider the weighted Bergman spaces
$$
\ahdp=\bigl\{f\in\hol(\D):\|f\|^p_{p,h}= \int_\D
|f(z)|^pe^{-ph(z)}dm_2(z)<\infty\bigr\},
$$
where $dm_2$ is area measure, $1\le p<\infty$.

A subset $\Gamma$ of $\D$ is a sampling set for $\ahdp$ if
$$
\|f\|^p_{p,h}\asymp\|f\|^p_{p,h,\Gamma}=\sum_{z\in\Gamma}e^{-ph(z)}
\frac{|f(z)|^p}{\Delta h(z)}.
$$

A subset $\Gamma$ of $\D$ is an interpolation set
for $\ahdp$ if there exists $\delta>0$ such that
for every function $a$ defined on $\Gamma$ such that $\|a\|_{p,h,\Gamma}<\infty$
there exists $f\in \ahdp$ such that
$$
a=f\bigm |  \Gamma, \qquad \delta \|f\|_{p,h}\le \|a\|_{p,h,\Gamma}.
$$

For a motivation of these definitions 
let us consider the case $p=2$. Then $\ahdd$ is a Hilbert 
space of analytic functions
in the unit disc, and we define the reproducing kernel
$k_\lambda\in\ahdd$, $\lambda\in\D$, by
$$
\langle k_\lambda,f\rangle=f(\lambda), \qquad f\in\ahdd.
$$
For regular $h$ considered in our paper (for precise
conditions on $h$ see the next section) we have
(see, for example, Lemmas~\ref{m10} and \ref{m1} (ii) below)
$$
\|k_\lambda\|^2\asymp e^{2h(\lambda)}\Delta h(\lambda).
$$
Therefore, a family of normalized reproducing kernels 
$\{k_\lambda/\|k_\lambda\|\}_{\lambda\in\Lambda}$
is a {\it frame} in $\ahdd$, that is
$$
\sum_{\lambda\in\Lambda}\Bigl|\bigl\langle f,
\frac{k_\lambda}{\|k_\lambda\|}\bigr\rangle\Bigr|^2\asymp
\|f\|^2,\qquad f\in\ahdd,
$$
if and only if $\Lambda$ is a sampling set for $\ahdd$.
In a similar way, a family $\{k_\lambda/\|k_\lambda\|\}_{\lambda\in\Lambda}$
is a {\it Riesz basic sequence} in its closed linear span in
$\ahdd$, that is
$$
\bigl\|\sum_{\lambda\in\Lambda}
a_\lambda 
\frac{k_\lambda}{\|k_\lambda\|} \Bigr\|^2\asymp
\sum_{\lambda\in\Lambda}
|a_\lambda|^2 
$$
for any sequence $\{a_\lambda\}_{\lambda\in\Lambda}$ of complex numbers,
if and only if $\Lambda$ is an interpolation set for $\ahdd$.

The famous Feichtinger conjecture (see, for example, \cite{CA1,CA2}) claims that any frame in a Hilbert space is a finite union
of Riesz basic sequences. For families of normalized reproducing kernels in $\ahdd$ this conjectrure 
translates into the question on whether
any sampling set for $\ahdd$ is a finite union of interpolation
sets for $\ahdd$. The answer is positive for $h$ we consider in this paper (as follows from Therems~\ref{tt01b} and 
\ref{tti01b}).
\medskip

In the plane case, if $h:[0,\infty)\to[0,+\infty)$ is an increasing
function such that $h(0)=0$, $\lim_{r\to\infty}h(r)=+\infty$, we
extend $h$ by $h(z)=h(|z|)$, $z\in\CC$, and consider the Banach
space $\ahc$ of entire functions with the norm
$$
\|f\|_h=\sup_{z\in\CC}|f(z)|e^{-h(z)}<+\infty,
$$
and the weighted Fock spaces
$$
\ahcp=\bigl\{f\in\hol(\CC):\|f\|^p_{p,h}=\int_\CC
|f(z)|^pe^{-ph(z)}dm_2(z)<\infty\bigr\}.
$$
and define the sampling and the interpolation subsets for the spaces $\ahc$, $\ahcp$,
$1\le p<\infty$, 
like above, in the disc case.
\medskip

K.~Seip and R.~Wallst\'en \cite{S,SW} described sampling and interpolation sets 
for the Fock spaces $\ahc$, $\ahcd$, with $h(z)=c|z|^2$, in terms of Beurling type densities. 
Later on, K.~Seip \cite{S1}
obtained such a description for the Bergman type spaces $\ahd$, 
$h(z)=\alpha\log\frac{1}{1-|z|}$, $\alpha>0$, and for $\mathcal{A}^2_0(\D)$ 
($=\ahdd$ with $h=0$).
For motivation and some applications of these results, for example to Gabor wavelets, 
see a survey \cite{BR} by J.~Bruna. 

The results of K.~Seip were extended to the Fock spaces $\mathcal{A}^p_h(\CC)$
(with $h$ not necessarily radial) such that $\Delta h\asymp 1$ in \cite{BO} and \cite{OS},
and to Bergman spaces $\mathcal{A}^p_h(\D)$, $\Delta h(z)\asymp (1-|z|^2)^{-2}$, in \cite{S2}. 
Yu.~Lyubarskii and K.~Seip \cite{LS} obtained such results for the spaces $\ahc$, $\ahcd$,
with $h(z)=m(\arg z)|z|^2$, $m$ being a $2\pi$ periodic $2$-trigonometrically convex function.
For more results and references see the books \cite{HKZ} and \cite{SB}.

Recently, N.~Marco, X.~Massaneda and J.~Ortega-Cerd\`a \cite{MMO} 
described sampling and interpolation sets 
for the Fock spaces $\mathcal{A}^p_h(\CC)$ for a wide class of $h$ such that
$\Delta h$ is a doubling measure. These results rely mainly upon the method used by A.~Beurling \cite{B} 
in his work on band-limited 
functions and on H\"ormander-type weighted estimates for the $\bar \partial$ equation. Therefore, it is not clear whether they can be extended to weight functions $h$ having more than
polynomial growth at infinity. 

The aim of our work is to extend previous results to the case 
of radial $h$ of arbitrary (more than polynomial) growth. For this, we use the method proposed 
by Yu.~Lyubarskii and K.~Seip in \cite{LS}. First we produce peak functions 
with precise asymptotics. For example, for every $z\in\D$ 
we find $f_z\in\ahd$ such that
$$
|f_z(w)|\asymp e^{h(w)-|w-z|^2\Delta h(z)/4}
$$
in a special neighborhood of $z$. (For a different type of peak functions in $\ahd$ see \cite{CS}.)
These peak functions permit us 
then to reduce our problems to those in the standard Fock spaces $\mathcal{A}^p_h(\CC)$, $h(z)=|z|^2$.

The construction of peak functions in \cite{LS} is 
based on sharp approximation of $h$ by $\log|f|$, $f\in\hol(\D)$,
obtained in the work of  
Yu.~Lyu\-barskii and M.~Sodin \cite{LSO}.
Here we need a similar construction for radial $h$ of arbitrary (more than polynomial) growth.
This is done in a standard way: we atomise the measure $\Delta h(z)dm_2(z)$ and
obtain a discrete measure $\sum \delta_{z_n}$. Since our $h$ are radial, we try 
to get sufficiently symmetric sequence $\{z_n\}$:
$$
\{z_n\}=\cup_k\{s_ke^{2\pi i m/N_k}\},\qquad s_k\to 1,\, N_k\to\infty.
$$
For approximation of general $h$ see the paper \cite{LM}
by Yu.~Lyubarskii and E.~Malinnikova and the references there.

Our paper is organized as follows. The main results are formulated in Section~\ref{mr}.
We construct peak functions in Section~\ref{sp}. Technical lemmas on sampling
and interpolation sets are contained in Section~\ref{ss}.
In Section~\ref{sld} we obtain auxiliary results on asymptotic densities. The theorems
on sampling sets are proved in Section~\ref{ssa}, and the theorems on interpolation sets
are proved in Section~\ref{sin}. In Sections~\ref{sp}--\ref{sin} we deal with the disc case.
Some changes necessary to treat the plane case are discussed in Section~\ref{spl}.
Finally, in Section~\ref{sindex} we give an application of our results
to subspaces of $\ahdp$ invariant under multiplication by the independent variable.

We do not discuss here the following interesting fact: the 
families of interpolation (sampling) sets are not monotonic 
with respect to the weight function $h$.
Also, we leave open other questions related to our results, including whether our interpolation sets are just sets of {\it free interpolation}, that is (say, for the spaces $\ahcp$) the sets $\Lambda\subset\CC$ such that
$$
\ell^\infty(\Lambda)\cdot \ahcp{\bigm|}\Lambda=
\ahcp{\bigm|}\Lambda.
$$
We hope to return to these questions later on.

The authors are grateful to Yu.~Lyubarskii, N.~Nikolski, K.~Seip, M.~Sodin, and P.~Thomas for helpful discussions.

\section{Main results}
\label{mr}

From now on in the disc case we assume that the function
$$
\rho(r)=\bigl[(\Delta h)(r)\bigr]^{-1/2},\qquad 0\le r<1,
$$
decreases to $0$ near the point $1$, and
\begin{equation}
\rho'(r)\to 0,\qquad r\to 1.
\label{n5}
\end{equation}
Then for any $K>0$, for $r\in(0,1)$ sufficiently close to $1$,
we have $[r-K\rho(r),r+K\rho(r)]\subset(0,1)$, and
\begin{equation}
\rho(r+x)=(1+o(1))\rho(r),\qquad |x|\le K\rho(r),\quad r\to 1.
\label{01}
\end{equation}

Furthermore, we assume that either \conod\ the function
$r\mapsto \rho(r)(1-r)^{-C}$
increases for some $C<\infty$ and for $r$ close to $1$ or 
\begin{equation*}
\text{\condd}\qquad\qquad \rho'(r)\log 1/\rho(r)\to 0,\qquad r\to 1.
\end{equation*}

Typical examples for \conod\ are
$$
h(r)=\log\log\frac 1{1-r}\cdot \log\frac 1{1-r}, \qquad h(r)=\frac 1{1-r},\qquad r\to1;
$$
a typical example for \condd\ is
$$
h(r)=\exp \frac 1{1-r},\qquad r\to1.
$$

Denote by $\mathcal D(z,r)$ the disc of radius $r$
centered at $z$, $\mathcal D(r)=\mathcal D(0,r)$.

Given $z,w\in\D$, we define
$$
d_\rho(z,w)=\frac{|z-w|}{\min(\rho(z),\rho(w))}.
$$
We say that a subset $\Gamma$ of $\D$ is $d_\rho$-separated (with constant $c$) if
$$
\inf\{d_\rho(z,w): z,w\in\Gamma,\,z\ne w\}\ge c>0.
$$

Given $\Gamma\subset\D$, we define its lower $d_\rho$-density
$$
\dmGD=\liminf_{R\to\infty}\liminf_{|z|\to 1,\,z\in\D}
\frac{\card(\Gamma\cap\mathcal D(z,R\rho(z)))}{R^2},
$$
and its upper $d_\rho$-density
$$
\dpGD=\limsup_{R\to\infty}\limsup_{|z|\to 1,\,z\in\D}
\frac{\card(\Gamma\cap\mathcal D(z,R\rho(z)))}{R^2}.
$$

We remark here that by \eqref{01}, for fixed $R$, $0<R<\infty$, we have
$$
\lim_{|z|\to 1,\,z\in\D}\frac1{\pi R^2}\int_{\mathcal D(z,R\rho(z))}\Delta h(z)\, dm_2(z)=1.
$$

We could compare these densities $D^\pm_{\rho}$ to those defined in the case $\rho(z)\asymp 1-|z|$
by K.~Seip in \cite{S1}:
\begin{align}
D^-(\Gamma,\D)&=\liminf_{r\to 1-}\inf_{\Phi\in\Aut(\D)}\frac{\sum_{z\in\Phi(\Gamma)\cap 
\mathcal D(r)\setminus\mathcal D(1/2)}\log\frac 1{|z| }}{\log\frac 1{1-r }},\\
D^+(\Gamma,\D)&=\limsup_{r\to 1-}\sup_{\Phi\in\Aut(\D)}\frac{\sum_{z\in\Phi(\Gamma)\cap 
\mathcal D(r)\setminus\mathcal D(1/2)}\log\frac 1{|z| }}{\log\frac 1{1-r }},
\end{align}
where $\Aut(\D)$ is the group of the M\"obius authomorphisms of the unit disc.
In contrast to $D^\pm$, the densities $D^\pm_{\rho}$ are rather ``local'', 
and correspondingly,
it is not difficult to compute $D^\pm_{\rho}(\Gamma,\D)$ for many concrete $\Gamma$.

\begin{thm} A set $\Gamma\subset\D$ is a sampling set for $\ahd$
if and only if it contains a $d_\rho$-separated subset
$\Gamma^{*}$ such that $D^{-}_\rho(\Gamma^{*},\D)> \frac 12$.
\label{tt01}
\end{thm}

\begin{thm} A set $\Gamma\subset\D$ is a sampling set for $\ahdp$, $1\le p<\infty$,
if and only if {\rm (i)} $\Gamma$ is a finite union of
$d_\rho$-separated subsets, and {\rm (ii)}
$\Gamma$ contains a $d_\rho$-separated subset $\Gamma^{*}$
such that
$D^{-}_\rho(\Gamma^{*},\D)> \frac 12$.
\label{tt01b}
\end{thm}

\begin{thm} A set  $\Gamma$ is an interpolation set for $\ahd$
if and only if it is $d_\rho$-separated and
$D^{+}_\rho(\Gamma,\D)< \frac 12$.
\label{tti01}
\end{thm}

\begin{thm} A set $\Gamma$ is an interpolation set for $\ahdp$, $1\le p<\infty$,
if and only if it is $d_\rho$-separated and
$D^{+}_\rho(\Gamma,\D)< \frac 12$.
\label{tti01b}
\end{thm}

In the plane case we assume that $h\in C^2(\CC)$, $\Delta h(z)\ge
1$, $z\in \CC$, that the function
$$
\rho(r)=\bigl[(\Delta h)(r)\bigr]^{-1/2},\qquad 0\le r<\infty,
$$
decreases to $0$ at infinity, and that
$$
\rho'(r)\to 0,\qquad r\to \infty.
$$
Then for any $K>0$, for sufficiently large $r$,
we have $r>K\rho(r)$, and
$$
\rho(r+x)=(1+o(1))\rho(r),\qquad |x|\le K\rho(r),\quad r\to \infty.
$$
Furthermore, we assume that
either \conoc\ the function
$r\mapsto \rho(r)r^{C}$
increases for some $C<\infty$ and for large $r$ or \condc\ 
$\rho'(r)\log 1/\rho(r)\to 0$ as $r\to\infty$.

Typical examples for \conoc\ are
$$
h(r)=r^2\log\log r, \qquad h(r)=r^{p},\,\,p>2,\qquad r\to \infty;
$$
a typical example for \condc\ is
$$
h(r)=\exp r,\qquad r\to \infty.
$$

Next, we introduce $d_\rho$, and the notion of $d_\rho$-separated subsets of $\CC$
as above. Given $\Gamma\subset\CC$, we define its lower $d_\rho$-density
$$
\dmGC=\liminf_{R\to\infty}\liminf_{|z|\to \infty}
\frac{\card(\Gamma\cap\mathcal D(z,R\rho(z)))}{R^2},
$$
and its upper $d_\rho$-density
$$
\dpGC=\limsup_{R\to\infty}\limsup_{|z|\to \infty}
\frac{\card(\Gamma\cap\mathcal D(z,R\rho(z)))}{R^2}.
$$

\begin{thm} A set $\Gamma\subset\CC$ is a sampling set for $\ahc$
if and only if it contains a $d_\rho$-separated subset $\Gamma^{*}$
such that $D^{-}_\rho(\Gamma^{*},\CC)> \frac 12$.
\label{tt02}
\end{thm}

\begin{thm} A set $\Gamma\subset\CC$ is a sampling set for $\ahcp$, $1\le p<\infty$,
if and only if {\rm (i)} $\Gamma$ is a finite union of
$d_\rho$-separated subsets and {\rm (ii)}
$\Gamma$ contains a $d_\rho$-separated subset $\Gamma^{*}$
such that $D^{-}_\rho(\Gamma^{*},\CC)> \frac 12$.
\label{tt02b}
\end{thm}

\begin{thm} A set  $\Gamma$ is an interpolation set for $\ahc$
if and only if it is $d_\rho$-separated and
$D^{+}_\rho(\Gamma,\D)< \frac 12$.
\label{tti02}
\end{thm}

\begin{thm} A set $\Gamma$ is an interpolation set for $\ahcp$, $1\le p<\infty$,
if and only if it is $d_\rho$-separated and
$D^{+}_\rho(\Gamma,\D)< \frac 12$.
\label{tti02b}
\end{thm}

\section{Peak functions}
\label{sp}

In this section we first approximate $h$ by $\log|f|$, for a special infinite 
product $f$. Then, using this construction, and an estimate on 
the partial products for the Weierstrass $\sigma$-function, we approximate 
the function $w\mapsto h(w)-|w-z|^2\Delta h(z)/4$
in a fixed $d_\rho$-neighborhood of $z\in\D$.

\begin{prop} There exist sequences $\{r_k\}$, $\{s_k\}$, $0=r_0<s_0<r_1<\ldots
r_{k}<s_{k}<r_{k+1}<\ldots <1$, and a sequence $N_k$, $k\ge 0$, of
natural numbers, such that $N_{k+1}\ge N_{k}$ for large $k$, and
\begin{itemize}
\item[(i)] \quad $\displaystyle \lim_{k\to\infty}\frac{r_{k+1}-r_k}{\rho(r_k)}=\sqrt{2\pi}$,
\quad  $\displaystyle \lim_{k\to\infty}N_k(r_{k+1}-r_k)=2\pi$,\newline\noindent
\phantom{A} $\displaystyle \lim_{k\to\infty}\frac {r_{k+1}-r_k}{r_{k}-r_{k-1}}=1$,
 \quad $\displaystyle \lim_{k\to\infty}\frac {r_{k+1}-s_k}{r_{k+1}-r_{k}}=\frac12$,
\item[(ii)] \quad if $\Lambda=\bigl\{s_ke^{2\pi i m/N_k}\bigr\}_{k\ge 0,\, 0\le m<N_k}$, and if
\begin{equation*}
f(z)=\lim_{r\to 1-}
\prod_{\lambda\in \Lambda\cap r\D}\Bigl(\frac{1-z/\lambda}{1-z\overline{\lambda}}\Bigr),
\end{equation*}
then the products in the right hand side converge uniformly on compact subsets of the unit disc, and
\begin{equation}
|f(z)|\asymp e^{h(z)}\frac{\dist(z,\Lambda)}{\rho(z)},\qquad z\in\D.\label{03}
\end{equation}
\item[(iii)] \quad Given $s\in(0,1)$ sufficiently close to the point $1$,
we can define $\{r_k\}$, $\{N_k\}$, $\{s_k\}$ and $\Lambda$ as above in such a way that
$s\in\Lambda$ and \eqref{03} holds {\rm(}uniformly in $s${\rm)}.
\end{itemize}
\label{ll6}
\end{prop}

\begin{proof} (i) We choose the sequence $\{r_k\}$ in the following way: $r_0=0$;
given $r_k$, $k\ge 0$, the number $r^*_k$ is defined by
$$
(r^*_k-r_k)\int_{r_k\le|z|<r^*_{k}}\Delta h(w)\frac{dm_2(w)}{2\pi}=2\pi,
$$
and $r_{k+1}$ is the smallest number in the interval $[r^*_k,1)$ such that
$$
N_k=\int_{r_k\le|z|<r_{k+1}} \Delta h(w)\frac{dm_2(w)}{2\pi}\in\NN.
$$
Here we use that by \eqref{n5},
\begin{equation}
\Delta h(r)^{-1/2}=\rho(r)=o(1-r), \qquad r\to 1.
\label{f84}
\end{equation}

Furthermore, since $\Delta h(r)$ increases for $r$ close to $1$, we obtain that
$N_k$ do not decrease for large $k$, and
$$
\lim_{k\to\infty}N_k(r_{k+1}-r_k)=2\pi.
$$
By \eqref{01} we obtain
\begin{gather*}
\lim_{k\to\infty}\frac{(r_{k+1}-r_k)^2}{\rho(r_k)^2}=2\pi,\\
\lim_{k\to\infty}\frac {r_{k+1}-r_{k}}{r_{k}-r_{k-1}}=1.
\end{gather*}

Next we define $s_k$ by the relations
\begin{equation}
\log\frac{1}{s_{k}}
=\frac 1{N_k} \int_{r_k\le |z|<r_{k+1}}\!\Delta h(w)\log\frac{1}{|w|}\frac{dm_2(w)}{2\pi},
\qquad k\ge 0.
\label{n8}
\end{equation}
Clearly,
$$
\log\frac{1}{r_{k+1}}<\log\frac{1}{s_{k}}<\log \frac{1}{r_{k}},\qquad k\ge 0.
$$
By \eqref{01} and \eqref{f84} we have
$$
\lim_{k\to\infty}\frac{r_{k+1}-s_k}{r_{k+1}-r_k}=\frac 12.
$$

(ii) First of all we note that 
$$
f(z)=\lim_{L\to \infty}\prod_{0\le m\le L} \frac{1-z^{N_m}s_m^{-N_m}}{1-z^{N_m}s_m^{N_m}}.
$$
Set 
$$
W_m=\log \Bigl|\frac{1-z^{N_m}s_m^{-N_m}}{1-z^{N_m}s_m^{N_m}}.
\Bigr|.
$$

If $m>k$, then for some constant $c>0$,
$$
s_m-r\ge c\frac{m-k}{N_m}.
$$
Since
\begin{equation}
\frac xy\le e^{x-y},\qquad 0\le x\le y\le 1,
\label{n14}
\end{equation}
we get
$$
r^{N_m}s_m^{-N_m}\le e^{-c(m-k)}.
$$

Next we use that if $|\zeta|\le c<1$, $|\zeta'|\le 1$, then
\begin{equation}
\Bigl|\log\Bigl|\frac{1-\zeta}{1-\zeta\zeta'}\Bigr|\Bigr|\le
c_1|\zeta|, \label{sp2}
\end{equation}
with $c_1$ depending only on $c$.

Therefore,
$$
|W_m|\le c_1 e^{-c(m-k)},\qquad m>k.
$$
Summing up, we obtain
\begin{equation}
\sum_{m>k}|W_m|\le c_1\sum_{m>k}e^{-c(m-k)}\le c_2,
\label{n4}
\end{equation}
for some positive constants $c_1,c_2$.

Thus, 
$$
f(z)=\prod_{m\ge 0} \frac{1-z^{N_m}s_m^{-N_m}}{1-z^{N_m}s_m^{N_m}}.
$$

Suppose that $z\in\D\setminus\Lambda$, $r_k\le r<r_{k+1}$, where $r=|z|$,
and for some $d$, $0\le d<N_k$,
\begin{equation*}
\bigl|\arg z-\frac{2\pi d}{N_k}\bigr|\le \frac{\pi}{N_k}.
\end{equation*}
Now we set
$$
A(z)=\log |f(z)|-h(z)-\log\frac{\dist(z,\Lambda)}{\rho(z)}.
$$

By Green's formula,
\begin{equation}
h(r)=\int_{\mathcal D(r)} \Delta h(w)\log\frac{r}{|w|}\frac{dm_2(w)}{2\pi}.
\label{n27}
\end{equation}
Therefore,
\begin{multline*}
A(z)+\log\frac{\dist(z,\Lambda)}{\rho(z)}=\\
\sum_{0\le m\le k-1}\biggl[ \log \Bigl|\frac {1-z^{N_m}s_m^{-N_m}}
{1-z^{N_m}s_m^{N_m}}\Bigr|-\int_{r_{m}\le |w|<r_{m+1}}\Delta
h(w)\log\frac{r}{|w|}\frac{dm_2(w)}{2\pi}\biggr]\\
-\int_{r_{k}\le |w|<r}\Delta h(w)\log\frac{r}{|w|}\frac{dm_2(w)}{2\pi}+
\sum_{m\ge k}\log \Bigl|\frac {1-z^{N_m}s_m^{-N_m}}
{1-z^{N_m}s_m^{N_m}}\Bigr|\\
=\sum_{0\le m\le k-1}U_m-V+\sum_{m\ge k} W_m.
\end{multline*}
\smallskip

First,
\begin{multline}
|V|\le \Bigl(\sup_{r_{k}\le |w|<r}\log\frac{r}{|w|}\Bigr) \int_{r_{k}
\le |w|<r}\Delta h(w)\frac{dm_2(w)}{2\pi}\\
\le cN_k(r_{k+1}-r_k)\le c_1,
\label{n3}
\end{multline}
for some constants $c,c_1$.
\smallskip

Next we are to verify that
\begin{equation}
\biggl|W_k-\log\Bigl|\frac{\dist(z,\Lambda)}{\rho(z)}\Bigr|\biggr|\le c,
\label{n2}
\end{equation}
for some constant $c$.
Indeed, $|z^{N_k}s_k^{N_k}|\le c_1$, for some constant $c_1<1$,
and it remains to estimate $W^*(z)$, where
$$
W^*(\zeta)=\log\Bigl|\frac{\rho(r_k)(1-\zeta^{N_k}s_k^{-N_k})}{\zeta-s_ke^{2\pi i d/N_k}}\Bigr|,
\qquad \zeta\in\D.
$$
Consider the set
$$
\Omega=\Bigl\{re^{i\theta}:r_k\le r\le
r_{k+1},\,\bigl|\theta-\frac{2\pi d}{N_k}\bigr|\le\frac{\pi}{N_k}\Bigr\}.
$$
For $\zeta\in\partial\Omega$ we have $\rho(r_k)\asymp
|\zeta-s_ke^{2\pi i d/N_k}|$, $|\zeta^{N_k}s_k^{-N_k}|\asymp 1$,
and either $|\arg(\zeta^{N_k}s_k^{-N_k})|\ge c_1$ or
$|1-|\zeta^{N_k}s_k^{-N_k}||\ge c_1$, for some positive
constant $c_1$. Therefore, $|W^*(\zeta)|\le c_2$, $\zeta\in
\partial\Omega$, for some constant $c_2$. Since $W^*$ is
harmonic on $\Omega$, we obtain, by the maximum
principle, that $|W^*(z)|\le c_2$, and \eqref{n2} follows. 
\smallskip

It remains to verify that
\begin{equation}
\sum_{0\le m<k}|U_m|\le c,
\label{04}
\end{equation}
with $c$ independent of $r$.
This together with \eqref{n3}, \eqref{n2} and \eqref{n4} implies 
that $A$ is bounded uniformly in $z\in\D\setminus\Lambda$, and \eqref{03} follows.
\smallskip

Since
$$
\int_{r_{m}\le |w|<r_{m+1}}\Delta h(w)\log\frac{r}{|w|}\frac{dm_2(w)}{2\pi}
=N_m\log\frac r{s_m},
$$
we have
$$ 
U_m=
\log \Bigl|\frac {1-
s_m^{N_m}z^{-N_m}}{1-z^{N_m}s_m^{N_m}}\Bigr|,\qquad 0\le m<k.
$$ 

Next we consider two cases. If $\rho$
satisfies the property \conod, then we define
$$
U^*_m(w)=\log \Bigl|\frac {1-s_m^{N_m}w^{-N_m}}{1-w^{N_m}s_m^{N_m}}\Bigr|,\qquad 0\le m<k,
\quad s_m<|w|\le 1,
$$
and divide $m$, $0\le m<k$, into the groups
$$
S_t=\bigl\{m:1-s_m\in\bigl[2^t(1-r),2^{t+1}(1-r)\bigr)\bigr\},\qquad t\in\Z_+.
$$
Put $\rho_t=\rho(1-2^t(1-r))$, $t\ge 0$.
Then by \conod\ we have
$$ 
\rho(s_m)\asymp \rho_t, \qquad m\in S_t,\quad t\in\Z_+,
$$
and hence,
$$
s_{m+1}-s_{m}\asymp \rho_t, \qquad m\in S_t,\quad t\in\Z_+.
$$ 
By \eqref{f84}, $\rho_t=o(2^t(1-r))$ as $2^t(1-r)\to 0$. Therefore,
for $t\ge 1$, $m\in S_t$, and for some $c<1$ independent of $m,r$,
we have by \eqref{n14} that
\begin{equation}
s_m^{N_m}(1-2^{t-1}(1-r))^{-N_m}\le c.
\label{sp1}
\end{equation}

Now, \eqref{sp2} implies that for
$\zeta\in(1-2^{t-1}(1-r))\T$,
\begin{gather*}
\bigl|U^*_m(\zeta)\bigr|\le c\Bigl(\frac{s_m}{1-2^{t-1}(1-r)}\Bigr)^{N_m}\\
\le c_1\exp\Bigl[-\frac{c_2((1-s_m)-2^{t-1}(1-r))}{\rho_t}\Bigr]
\le c_1\exp\Bigl[-\frac{c_2 2^{t-1}(1-r)}{\rho_t}\Bigr],
\end{gather*}
and
\begin{equation}
\bigl|\sum_{m\in S_t}U^*_m(\zeta)\bigr|\le
c_3\frac{2^t(1-r)}{\rho_t}\exp\Bigl[-\frac{c_22^{t-1}(1-r)}{\rho_t}\Bigr]\le
c_4,\qquad t\ge 1,
\label{f29}
\end{equation}
with $c,c_1,c_2,c_3,c_4$ independent of $r$, $t$.

Furthermore,
$$
\bigl|\sum_{m\in S_0}U^*_m(r)\bigr|\le 
c\sum_{m\in S_0}\Bigl(\frac{s_m}{r}\Bigr)^{N_m}\le 
c\sum_{m\in S_0}e^{-c_1(k-m)}\le c_2,
$$
with positive $c,c_1,c_2$ independent of $r$.

Since $U^*_m(w)=0$, $w\in\T$, and for $t\ge 1$, $m\in S_t$, $U^*_m$ are harmonic in the annulus
$\{w: 1-2^{t-1}(1-r)\le w\le 1\}$, we deduce from \eqref{f29} that
$$
\bigl|\sum_{m\in S_t}U^*_m(r)\bigr|\le c2^{-t},\qquad t\ge 1,
$$
with $c$ independent of $r,t$, and \eqref{04} is proved.
\smallskip

Suppose now that $\rho$ satisfies the property \condd. For some constant $c>0$,
$$
r-s_m\ge \frac c{N_m},\qquad m<k,
$$
and we obtain that
$$
(s_m/r)^{N_m}\le e^{-c},
$$
and again by \eqref{sp2},
$$
|U_m|\le c_1(s_m/r)^{N_m},\qquad m<k,
$$
for some positive constant $c_1$.

By \eqref{01} and \eqref{n14} we obtain that
\begin{multline}
\sum_{m<k}|U_m|\le
c\sum_{m<k}e^{N_m(s_m-r)}\\  \le
c_1\sum_{m<k}\int_{r_m}^{r_{m+1}}e^{-c_2(r-x)/\rho(x)}\frac{dx}{\rho(x)}\le
c_1\int_0^re^{-c_2(r-x)/\rho(x)}\frac{dx}{\rho(x)},
\label{05}
\end{multline}
with positive $c,c_1,c_2$ independent of $k$, $r$.

Choose $y$ such that $\rho(y)=2\rho(r)$. Then
for $x<y$ close to $1$,
\begin{gather*}
\rho(x)\log\frac1{\rho(x)}\le
\rho(r)\log\frac1{\rho(r)}+(r-x)\sup_{[x,r]}\bigl|\bigl[\rho\log\frac1{\rho}\bigr]'\bigr|,\\
\rho(x)\log\frac1{\rho(x)}\le c_2(r-x).
\end{gather*}
Hence,
\begin{equation}
\int_0^y e^{-c_2(r-x)/\rho(x)}\frac{dx}{\rho(x)}\le 
c_3+\int_0^y e^{-\log(1/\rho(x))}\frac{dx}{\rho(x)}\le c_3+1,
\label{h4}
\end{equation}
with $c_3$ independent of $r$.

Finally, for $r$ close to $1$
$$
\int_y^r e^{-c_2(r-x)/\rho(x)}\frac{dx}{\rho(x)}\le
\int_y^re^{-c_2(r-x)/(2\rho(r))}\frac{dx}{\rho(r)}\le c_3,
$$
with $c_3$ independent of $r$.
These inequalities together with \eqref{05} prove \eqref{04}, and hence, \eqref{03}.

(iii) Given $s_k\le s< s_{k+1}$, $k\ge 0$, we may find $0<r'_k<s<r'_{k+1}<1$
such that
$$
N_k=N'_k=\int_{r'_k\le|z|<r'_{k+1}} \Delta
h(w)\frac{dm_2(w)}{2\pi},
$$
and
$$
\log\frac{1}{s}
=\frac 1{N_k} \int_{r'_k\le |z|<r'_{k+1}}\!\Delta h(w)\log\frac{1}{|w|}\frac{dm_2(w)}{2\pi}.
$$
After that, we define $r^{\prime *}_n$, $r'_n$, $N'_{n-1}$,
$n>k+1$, as in part (i).

Furthermore, we define by induction, on the step $t\ge 1$,
the number $r^{\prime *}_{k-t+1}\in (0, r^{\prime}_{k-t+1})$ by the equality
$$
(r^{\prime}_{k-t+1}-r^{\prime *}_{k-t+1})
\int_{r^{\prime *}_{k-t+1}\le|z|<r^{\prime}_{k-t+1}}\!\!\!\Delta h(w)\frac{dm_2(w)}{2\pi}=2\pi,
$$
and the number $r^{\prime}_{k-t}$ as the largest number in the interval $(0,r^{\prime *}_{k-t+1}]$
such that
$$
N'_{k-t}=\int_{r'_{k-t}\le|z|<r'_{k-t+1}} \!\!\!\Delta h(w)\frac{dm_2(w)}{2\pi}\in\NN.
$$

We continue this induction process until either
$$
A_p=r'_p\int_{|z|<r'_{p}} \Delta h(w)\frac{dm_2(w)}{2\pi}< 2\pi
$$
or $A_p\ge 2\pi$ and
$$
\int_{r<|z|<r'_{p}} \Delta h(w)\frac{dm_2(w)}{2\pi}\notin\NN, \qquad r\in [0,r^{\prime *}_{p}].
$$
It is clear that in both cases $r'_p\le c(h)<1$.
Next, we modify $h$ on $r'_{p+1}\D$ in such a way that the modified function
$h^*$ is smooth, radial, subharmonic, $|h^*(0)|\le c_1(h)$, and
$$
N=\int_{|z|<r'_{p+1}} \Delta h^*(w)\frac{dm_2(w)}{2\pi}\in\NN.
$$

Finally, we set $r_0=0$, $N_0=N$, $r_m=r'_{m+p}$, $N_m=N'_{m+p}$,
$m\ge 1$, and define $s_m$, $m\ge 0$, by \eqref{n8}. We apply the
above argument to $h=h^*-h^*(0)$ to obtain all the estimates from
(i)--(ii) uniformly in $s$ together with the property
$s\in\Lambda$.
\end{proof}

Given $0<r\le 1$, we define
\begin{equation}
h_r(w)=
\begin{cases}
h(w), & |w|< r,\\
h(r)+\log\frac{|w|}{r}\int_{\mathcal D(r)}\Delta h(z)\frac{dm_2(z)}{2\pi}, & |w|\ge  r.
\end{cases}
\label{n54}
\end{equation}
Note that
$$
h_r(w)=\int_{\mathcal D(r)}\Delta h(z)\bigl(\log\frac{|w|}{|z|}\bigr)\frac{dm_2(z)}{2\pi}, \qquad|w|\ge  r.
$$

The proof of Proposition~{\rm\ref{ll6}} gives us immediately

\begin{lem} In the notations of Proposition~{\rm\ref{ll6}}, if $s_{k-1}\le r<s_k$, and
$$
f(\zeta)=\prod_{\lambda\in \Lambda\cap \mathcal D(r)}
\Bigl(\frac{1-\zeta/\lambda}{1-\zeta\overline{\lambda}}\Bigr),
$$
then
$$
|f(\zeta)|\asymp e^{h_{r_k}(\zeta)}
\min\Bigl(1,\frac{\dist(\zeta,\Lambda\cap r\D)}{\rho(\zeta)}\Bigr),\qquad \zeta\in\D.
$$
\label{m8}
\end{lem}

If $|z|=r$, then we can divide this $f$ by three factors
$\zeta-\lambda_j$, $\lambda_j\in\Lambda\cap r\D\cap \mathcal
D(z,5\rho(z))$, $j=1,2,3$, and multiply it by $\rho(z)^3$, to obtain

\begin{lem} Given $z\in\D$ such that $r=|z|$ is sufficiently close to $1$, there
exists a function $g_z$ analytic and bounded in $\D$ and such that
uniformly in $z$,
\begin{align}
|g_z(w)|e^{-h(w)}&\asymp 1,\qquad |w-z|<\rho(z),\label{n30}\\
|g_z(w)|e^{-h(w)}&\le
c(h)\min\Bigl[1,\frac{\min[\rho(z),\rho(w)]}{|z-w|}\Bigr]^3,\qquad
w\in\D.\label{n31}
\end{align}
\label{m10}
\end{lem}

We need only to verify that for some $c$,
$$
e^{h_r(w)-h(w)}\rho(z)^3\le c\rho(w)^3, \qquad 0\le r=|z|\le |w|<1.
$$
This follows from the inequality $h(t)\ge h_r(t)$ and the estimate
$$
\frac{d}{dt}[h(t)-h_r(t)]=\frac{1}{t}\int_{\mathcal
D(t)\setminus\mathcal D(r)}\Delta h(w)\frac{dm_2(w)}{2\pi}\ge
\frac{3}{\rho(t)},\qquad t\ge r+B\rho(r),
$$
for some $B>0$ independent of $r$.

Next, we obtain an asymptotic estimate for partial products of the
Weierstrass $\sigma$-function.

\begin{lem} Given $R\ge 10$, we define $\Sigma=\Sigma_R=(\Z+i\Z)\cap \mathcal D(R^2)$,
$$
P_R(z)=z\prod_{\lambda\in
\Sigma\setminus\{0\}}\Bigl(1-\frac{z}{\lambda}\Bigr)
$$
Then uniformly in $R$
\begin{align}
|P_R(z)|&\asymp \dist(z,\Sigma)e^{(\pi/2)|z|^2},\qquad |z|\le R,\label{07}\\
|P_R(z)|&\ge
c\dist(z,\Sigma)\Bigl(\frac{e|z|^2}{R^2}\Bigr)^{(\pi/2)R^2},\qquad
|z|> R.\label{08}
\end{align}
\label{ll8}
\end{lem}

\begin{proof}
For every $\lambda\in \Sigma$ denote
$$
Q_\lambda=\{w\in\CC:|\RE(w-\lambda)|<1/2,\,|\IM(w-\lambda)|<1/2\}.
$$
Set
$$
Q=\bigsqcup_{\lambda\in \Sigma}Q_\lambda.
$$
If $z\in Q$, then we denote by $\lambda_0$ the element of $\Sigma$
such that $z\in Q_{\lambda_0}$. For
$\lambda\in\Sigma\setminus\{0\}$ we define
\begin{multline*}
B_\lambda=\int_{Q_\lambda}\log\Bigl|1-\frac zw\Bigr|\,dm_2(w)-\log\Bigl|1-\frac z\lambda\Bigr|\\
=\int_{Q_0}\Bigl[\log\Bigl|\frac {\lambda-z+w}{\lambda-z}\Bigr|
-\log\Bigl|\frac {\lambda+w}{\lambda}\Bigr|\Bigr]\,dm_2(w).
\end{multline*}
We use that
$$
\bigl|\log|1+a|-\RE (a-\frac{a^2}{2}) \bigr|=O(|a|^3),\qquad a\to 0.
$$
Since
$$
\int_{Q_0}w\,dm_2(w)=0,\qquad
\int_{Q_0}w^2\,dm_2(w)=0,
$$
we conclude that
$$
|B_\lambda|\le c\Bigl[\frac1{|\lambda|^3}+\frac1{|z-\lambda|^3}\Bigr],
\qquad \lambda\in\Sigma\setminus\{0,\lambda_0\}.
$$
Furthermore, we define
\begin{gather*}
B_0=\int_{Q_0}\log\Bigl|1-\frac zw\Bigr|\,dm_2(w)-\log |z|\\=
\int_{Q_0}\Bigl[\log\frac{1}{|w|}+\log\Bigl|1-\frac wz\Bigr|\Bigr]\,dm_2(w).
\end{gather*}
If $\lambda_0\ne 0$,
then $|B_0|\le c$ for an absolute constant $c$. Similarly, in this case,
\begin{gather*}
B_{\lambda_0}+\log|z-\lambda_0|=
\int_{Q_0}\Bigl[\log|w+(\lambda_0-z)|-\log\Bigl|1+\frac w{\lambda_0}\Bigr|\Bigr]\,dm_2(w),
\end{gather*}
and, hence, $\bigl|B_{\lambda_0}+\log|z-\lambda_0|\bigr|\le c$ for
an absolute constant $c$. In the same way, if $\lambda_0=0$, then
$\bigl|B_{0}+\log|z|\bigr|\le c$ for an absolute constant $c$.

Therefore,
\begin{equation}
\left.\!\!\!\!\!\!
\begin{aligned}
\Bigl|
\int_{Q}\log\Bigl|1-\frac zw\Bigr|\,dm_2(w)-
\log \frac{|P_R(z)|}{\dist(z,\Sigma)}
\Bigr|&=O(1),
\, |z|<R^2+1,\\
\Bigl|
\int_{Q}\log\Bigl|1-\frac zw\Bigr|\,dm_2(w)-
\log |P_R(z)|
\Bigr|&=O(1),
\, |z|\ge R^2+1.
\end{aligned}
\right\}
\label{n7}
\end{equation}

Next we use the identity
\begin{multline}
\frac{2}{\pi}\int_{\mathcal D(R^2+1)}\log\Bigl|1-\frac zw\Bigr|\,dm_2(w) 
=4\int_0^{\min(|z|,R^2+1)}\log\frac{|z|}{s}\,s\,ds
\\
\quad=\begin{cases}\displaystyle
(R^2+1)^2+2(R^2+1)^2\log\Bigl|\frac {z}{R^2+1}\Bigr|,\quad |z|\ge
R^2+1,
\\ |z|^2,\quad |z|< R^2+1,
\end{cases}
\label{09}
\end{multline}
and the estimates
\begin{multline}
\biggl|\int_{\mathcal D(R^2+1)\setminus Q}\log\Bigl|1-\frac zw\Bigr|\,dm_2(w)\biggr|\\
\le \biggl|\int_{\mathcal D(R^2+1)\setminus Q}\RE\frac zw\,dm_2(w)\biggr|+
c\int_{\mathcal D(R^2+1)\setminus Q}\Bigl|\frac zw\Bigr|^2dm_2(w)\\
\le \frac {c|z|^2}{R^2},\qquad |z|\le R^{3/2},\label{10}
\end{multline}
and
\begin{multline}
\biggl|\int_{\mathcal D(R^2+1)\setminus Q}\log\Bigl|1-\frac zw\Bigr|\,dm_2(w)\biggr|\\
\le m_2\bigl(\mathcal D(R^2+1)\setminus Q\bigr)\cdot \log\bigl[(R^2+1)(|z|+(R^2+1))\bigr]
,\,\,\, z\in\CC,\label{12}
\end{multline}
for an absolute constant $c$.
Now, \eqref{07} follows from \eqref{n7}--\eqref{10};
\eqref{08} follows from \eqref{n7}--\eqref{12}.
\end{proof}

\begin{prop} Given $R\ge 100$, there exists $\eta(R)>0$ such that for every $z\in\D$ with
$|z|\ge 1-\eta(R)$, there exists a function $g=g_{z,R}$ analytic
in $\D$ such that uniformly in $z,R$ we have
\begin{align}
|g(w)|e^{-h(w)}&\asymp e^{-|z-w|^2/[4\rho(z)^2]},\qquad
w\in\D\cap \mathcal D(z,R\rho(z)),\label{15}\\
|g(w)|e^{-h(w)}&\le
c(h)\Bigl[\frac{R^2\rho(z)^2}{e|z-w|^2}\Bigr]^{R^2/4},\qquad
w\in\D\setminus \mathcal D(z,R\rho(z)).\label{16}
\end{align}
\label{ll7}
\end{prop}

\begin{proof} Without loss of generality we may assume that $z\in (0,1)$.
By Proposition~\ref{ll6},
we find $\{r_k\}$, $\{N_k\}$, $\{s_k\}$, $\Lambda=\bigl\{s_ke^{2\pi i m/N_k}\bigr\}$,
and $f\in \ahd$ such that $z=s_k$ for some $k$, $Z(f)=\Lambda$, and
\begin{equation}
|f(w)|\asymp e^{h(w)}\frac{\dist(w,\Lambda)}{\rho(w)},\qquad w\in\D.
\label{n16}
\end{equation}
We define $\Sigma=\Sigma_{R/\sqrt{2\pi}}$ (see Lemma~\ref{ll8}),
and denote 
$$ 
\lambda_{a,b} =s_{k+a}e^{2\pi i b/N_{k+a}},\qquad a+bi\in\Sigma. 
$$
Then
$$
\max_{a+bi\in\Sigma}\Bigl|a+bi-\frac{\lambda_{a,b}-\lambda_{0,0}}
{\sqrt{2\pi}\rho(z)}\Bigr|\le \varepsilon(k),
$$
with $\varepsilon(k)\to 0$ as $k\to\infty$. Denote
$$
Q(w)=\frac{w-z}{\rho(z)}\prod_{a+bi\in\Sigma\setminus\{0\}}
\Bigl(\frac{w-\lambda_{a,b} } {z-\lambda_{a,b} }\Bigr),
$$
and put $g=f/Q$. Now, estimates \eqref{15} and \eqref{16} follow
for fixed $R$ when $k$ is sufficiently large, and correspondingly,
$\varepsilon(k)$ is sufficiently small.

Indeed, by \eqref{n16}, for $u=ge^{-h}$ we have
$$
|u(w)|\asymp
\frac{\dist(w,\Lambda)}{\rho(w)}\cdot\frac{\rho(z)}{|w-z|}\cdot
\prod_{a+bi\in\Sigma\setminus\{0\}}\Bigl|\frac
{z-\lambda_{a,b} }{w-\lambda_{a,b} }\Bigr|.
$$
If $w=z+\sqrt{2\pi}\rho(z)w'$, then
$$
|u(w)|\asymp \frac{\dist(w,\Lambda)}{\rho(w)}\cdot\frac 1{|w'|}
\cdot\!\prod_{a+bi\in\Sigma\setminus\{0\}}
\Bigl|\frac{(z-\lambda_{a,b} )/(\sqrt{2\pi}\rho(z))}
{w'+(z-\lambda_{a,b} )/(\sqrt{2\pi}\rho(z))}\Bigr|.
$$
For small $\varepsilon(k)$ and for $\dist(w',\Sigma)>1/10$ we have
$$
|u(w)|\asymp \frac{\dist(w,\Lambda)}{\rho(w)}\cdot\frac 1{|w'|}
\cdot\prod_{a+bi\in\Sigma\setminus\{0\}}
\Bigl|\frac{a+bi}{w'-(a+bi)}\Bigr|=
\frac{\dist(w,\Lambda)}{\rho(w) |P_{R/\sqrt{2\pi}}(w')|}.
$$
Now, for $z$ sufficiently close to $1$, by Lemma~\ref{ll8} and by the maximum principle,
we have
\begin{multline*}
|u(w)|\asymp
\frac{\dist(w',\Sigma)}{\bigl|P_{R/\sqrt{2\pi}}(w')\bigr|}\asymp
\exp\Bigl[-\frac{\pi}{2}|w'|^2\Bigr]\\
=\exp\Bigl[-\frac{|z-w|^2}{4\rho(z)^2}\Bigr],\qquad |w'|\le \frac {R}{\sqrt{2\pi}},
\end{multline*}
and
\begin{multline*}
|u(w)|\le
c\frac{\dist(w',\Sigma)}{\bigl|P_{R/\sqrt{2\pi}}(w')\bigr|}\le
c_1\Bigl[\frac{R^2}{2\pi e|w'|^2}\Bigr]^{(\pi/2)\cdot(R^2/(2\pi))}\\
=c_1\Bigl[\frac{R^2\rho(z)^2}{e|z-w|^2}\Bigr]^{R^2/4}, \qquad |w'|>\frac {R}{\sqrt{2\pi}}.
\end{multline*}

\end{proof}

\begin{prop} Given $R\ge 100$, there exists $\eta(R)>0$ such that for every $z\in\D$ with
$|z|\ge 1-\eta(R)$, there exists a function $g=g_{z,R}$ analytic
in $\D$ such that uniformly in $z,R$ we have
\begin{align*}
|g(w)|e^{-h(w)}&\asymp e^{-|z-w|^2/[4\rho(z)^2]},\qquad
w\in\D\cap \mathcal D(z,R\rho(z)),\\
|g(w)|e^{-h(w)}&\le
c(h)\Bigl[\frac{R^2\min[\rho(z),\rho(w)]^2}{e|z-w|^2}\Bigr]^{R^2/4},\quad
w\in\D\setminus \mathcal D(z,R\rho(z)).
\end{align*}
\label{m14}
\end{prop}

\begin{proof} We use the argument from the above proof,
and just replace Proposition~\ref{ll6} by Lemma~\ref{m8}.
Furthermore, we use the argument from the proof of
Lemma~\ref{m10}.
\end{proof}
\smallskip

\section{$d_\rho$-separated sets}
\label{ss}

Here we establish several elementary properties of $d_\rho$-separated sets,
sets of sampling, and sets of interpolation.

\begin{lem} Let $0<R<\infty$, let $z$ be sufficiently close to the unit circle, 
$\eta^*(R)<|z|<1$, and let $f$ be bounded and analytic in
$D=\mathcal D(z,R\rho(z))$. Then
\begin{multline*}
{\rm(i)}\quad \Bigl||f(z_1)|e^{-h(z_1)}-|f(z_2)|e^{-h(z_2)}\bigr|
\le c(R,h)\,d_\rho(z_1,z_2)\max_D|fe^{-h}|,\\
z_1,z_2\in \mathcal D(z,R\rho(z)/2),
\end{multline*}
\begin{gather*}
{\rm(ii)}\qquad\qquad |f(z)|e^{-h(z)}\le \frac{c(R,h)}{\rho(z)^2}\int_D |f(w)|e^{-h(w)}dm_2(w).
\end{gather*}
\label{m1}
\end{lem}

\begin{proof}We may assume that $\rho(\zeta)\asymp \rho(z)$, $\zeta\in D$.
We suppose that $\max_D|fe^{-h}|=1$ and define
$$
H(w)=h(z+wR\rho(z)),\qquad |w|\le 1.
$$
Then
$$
\Delta H(w)=\frac{R^2\rho(z)^2}{\rho(z+wR\rho(z))^2}\asymp R^2,\qquad |w|\le 1.
$$
Set
$$
G(w)=\int_\D\log\Bigl|\frac{z-w}{1-\bar zw}\Bigr|\Delta H(z)\,
dm_2(z),\qquad |w|\le 1.
$$
Then $|G(w)|+|\nabla G(w)|\le c$, $w\le 1$, for some $c$ depending
only on $h$ and $R$, and $H_1=H-G$ is real and harmonic in $\D$.
Denote by $\tilde H_1$ the harmonic conjugate of $H_1$, and
consider
$$
F(w)=f(z+wR\rho(z))e^{-H_1(w)-i\tilde H_1(w)}.
$$
Then $F$ is analytic and bounded in $\D$, and hence,
$$
|F(w_1)-F(w_2)|\le c|w_1-w_2|, \qquad w_1,w_2\in \mathcal D(1/2).
$$
Since
$$
|F(w)|=|f(z+wR\rho(z))|e^{-h(z+wR\rho(z))}e^{G(w)},
$$
we obtain assertion (i). Assertion (ii) follows by the mean value
property for $F$.
\end{proof}

\begin{cor} Every set of sampling for $\ahd$ contains a $d_\rho$-separa\-ted
set of sampling for $\ahd$.
\label{m2}
\end{cor}

\begin{cor} Every set of interpolation for $\ahd$ is $d_\rho$-separa\-ted.
\label{m15}
\end{cor}

\begin{cor} Every set of interpolation for $\ahdp$, $1\le p<\infty$, is $d_\rho$-separa\-ted.
\label{m16}
\end{cor}

\begin{lem} For every $\ve>0$, $1\le p<\infty$, we have 
$\ahdp\subset \mathcal{A}_{(1+\ve)h}(\D)$.
\label{m9}
\end{lem}

\begin{proof} By \eqref{n5} and \eqref{n27},
$$
\frac{|\rho'(r)|}{\rho(r)}=o\Bigl(\frac 1{\rho(r)}\Bigr)=o\bigl(h'(r)\bigr), \qquad r\to 1,
$$
and hence,
\begin{equation}
e^{\ve h(z)}\rho(z)^2\to\infty, \qquad |z|\to 1.
\label{j41}
\end{equation}
Applying H\"older's inequality and Lemma~\ref{m1}~(ii) with $R=1$, we obtain our assertion:
\begin{multline*}
|f(z)|e^{-(1+\ve)h(z)}\le\frac{c\cdot e^{-\ve h(z)}}{\rho(z)^{2}}
\int_{\mathcal D(z,\rho(z))}|f(z)|e^{-h(z)}dm_2(z)\\
\le
\frac{c\cdot e^{-\ve h(z)}}{\rho(z)^{2/p}}
\Bigl(\int_{\mathcal D(z,\rho(z))}|f(z)|^pe^{-ph(z)}dm_2(z)\Bigr)^{1/p}\le c,\qquad z\in\D.
\end{multline*}
\end{proof}

\begin{lem} Let $\Gamma$ be a $d_\rho$-separated {\rm(}with constant $\gamma${\rm)}
subset of $\D$. If $R>0$,  
$\Omega(R,\Gamma)=\{w:\min_{z\in\Gamma}d_\rho(w,z)\le R\}$, and if
$f$ is analytic in $\Omega(R,\Gamma)$, then
$$
\|f\|^p_{p,h,\Gamma}\le c(\gamma,R,h,p)\int_{\Omega(R,\Gamma)}|f(w)|^pe^{-ph(w)}dm_2(w).
$$
\label{m3}
\end{lem}

\begin{proof} The assertion follows from Lemma~\ref{m1}~(ii).
\end{proof}

\begin{lem} Let $\Gamma\subset\D$. Then
\begin{equation}
\|f\|^p_{p,h,\Gamma}\le c(\Gamma)\|f\|^p_{p,h},\qquad f\in\ahdp,
\label{n25}
\end{equation}
if and only if $\Gamma$ is a finite union of
$d_\rho$-separated subsets.
\label{m4}
\end{lem}

\begin{proof}
For every $z\in \D$ with $|z|$ close to $1$, 
we apply Lemma~\ref{m10} to obtain the function $f=g_{z}$
such that
\begin{align}
|f(w)|e^{-h(w)}&\asymp 1,\qquad |w-z|<\rho(z),\label{n37a}\\
|f(w)|e^{-h(w)}&\le
c(h)\frac{\rho(z)^3}{|z-w|^3},\qquad
w\in\D.\label{n37b}
\end{align}
By \eqref{n37a},
\begin{multline}
\|f\|^p_{p,h,\Gamma}\ge 
\sum_{w\in\Gamma \cap \mathcal D(z,\rho(z))}|f(w)|^pe^{-ph(w)}\rho(w)^2\\
\ge c\card \bigl(\Gamma\cap\mathcal D(z,\rho(z))\bigr)\rho(z)^2.\label{n37c}
\end{multline}
Furthermore, by \eqref{n37a}--\eqref{n37b},
\begin{align*}
\int_{|w-z|<\rho(z)}|f(w)|^pe^{-ph(w)}dm_2(w)&\asymp \rho(z)^2,\\
\int_{|w-z|\ge \rho(z)}|f(w)|^pe^{-ph(w)}dm_2(w)&\le c(h)\rho(z)^2,
\end{align*}
and hence, $f\in\ahdp$ and 
\begin{equation}
\|f\|^p_{p,h}\asymp \rho^2(z).
\label{n37d}
\end{equation}
Now, \eqref{n25}, \eqref{n37c}, and \eqref{n37d} imply that
$$
\sup_{z\in \D}\card \bigl(\Gamma\cap\mathcal D(z,\rho(z))\bigr) <\infty,
$$
and hence, $\Gamma$ is a finite union of
$d_\rho$-separated subsets.

In the opposite direction, if $\Gamma$ is
$d_\rho$-separated, then \eqref{n25} follows from Lemma~\ref{m3}.
\end{proof}

\begin{lem} Every set of sampling for $\ahdp$ contains a $d_\rho$-separa\-ted
set of sampling for $\ahdp$.
\label{m6}
\end{lem}

\begin{proof} Let $\Gamma$ be a set of sampling for $\ahdp$. For every $\ve>0$ we can find
a $d_\rho$-separa\-ted subset $\Gamma^*$ of $\Gamma$ such that
$$
\sup_{w\in\Gamma}\min_{z\in \Gamma^*} d_{\rho}(z,w)\le \ve.
$$
Suppose that there exists $f\in\ahdp$ such that
$$
\|f\|_{2,h}\asymp 1,\qquad \|f\|_{p,h,\Gamma}\asymp 1,\qquad \|f\|_{p,h,\Gamma^*}\le \ve.
$$
By Lemma~\ref{m4}, for some $N,K$ independent of $\ve$, both
$\Gamma$ and $\Gamma^*$ are unions of $N$ subsets
$d_\rho$-separated with constant $K$. 
Without loss of regularity we can assume that 
$|z|+\rho(z)<1$, $z\in \Gamma$.
For every $z_k\in\Gamma^*$
we choose $w_k\in \mathcal D(z_k,\rho(z_k))$ such that
$$
2|f(w_k)|^pe^{-ph(w_k)}\ge u_k^p=\sup_{\mathcal D(z_k,\rho(z_k))}|f|^pe^{-ph}.
$$
Then the sequence $\{w_k\}$ is the union of
$c(N,K)$ subsets $d_\rho$-separated with constant $c_1(N,K)$.
By Lemma~\ref{m3},
$$
\sum_{z_k\in\Gamma^*} u^p_k\rho(z_k)^2\le C\|f\|^p_{p,h},
$$
with $C$ independent of $\ve$. Furthermore, by Lemma~\ref{m1}~(i),
for every $k$ and for every $w\in \mathcal D(z_k,\ve\rho(z_k))$,
$$
\bigl||f(w)|e^{-h(w)}-|f(z_k)|e^{-h(z_k)}\bigr|\le C\ve u_k.
$$
Therefore,
\begin{multline*}
\|f\|^p_{p,h,\Gamma}=\sum_{w\in\Gamma}|f(w)|^pe^{-ph(w)}\rho(w)^2\\ \le
C\sum_{z_k\in\Gamma^*}\bigl(|f(z_k)|^pe^{-ph(z_k)}+\ve^pu^p_k\bigr) \rho(z_k)^2\\
\le
c(\|f\|^p_{p,h,\Gamma^*}+\ve^p\|f\|^p_{p,h})\le C\ve^p,
\end{multline*}
with $C$ independent of $\ve$. This contradiction implies our assertion.
\end{proof}

\section{Asymptotic densities}
\label{sld}

Given a set $\Gamma\subset\D$ such that $D^{-}_\rho(\Gamma,\D)<\infty$,
we denote
$$
q_{-}(R)=\liminf_{|z|\to 1,\, z\in\D}\frac{\card \Gamma(z,R)}{R^2}, \qquad 0<R<\infty,
$$
where
$$
\Gamma(z,R)=\Gamma\cap \mathcal D(z,R\rho(z)).
$$

In this section we study the behavior of the function $q_{-}$ and obtain 
a Beurling type result (Lem\-ma~\ref{ll9}).

We use the following

\begin{lem} \begin{itemize}
\item[(i)] If $R>0$, and $0<\varepsilon<\varepsilon(R)$, then
for $R''\ge R'(R,\varepsilon)$, and for $z\in\D$ such that $|z|\ge
\eta_1(R'')$ we have
$$
\frac{\card \Gamma(z,R'')}{R^{\prime\prime2}}\ge q_{-}(R)-\ve;
$$
\item[(ii)] if $\delta>0$, $R>0$, 
$R''\ge R'(R,\delta)$, $|z|\ge\eta_2(R'')$,
$$
E=\bigl\{w\in \mathcal D(z,R''\rho(z)/2): \frac{\card \Gamma(w,R)}{R^2}\ge q_{-}(R)+\delta\bigr\},
$$
and
$$
m_2E\ge \delta m_2\mathcal D(z,R''\rho(z)/2),
$$
then
$$
\frac{\card \Gamma(z,R'')}{R^{\prime\prime2}}\ge q_{-}(R)+\frac{\delta^2}{5}.
$$
\end{itemize}
\label{subl}
\end{lem}

\begin{proof} By \eqref{01},
$$
\max_{w\in \mathcal \mathcal D(z,R''\rho(z))}
\Bigl|\log\frac{\rho(z)}{\rho(w)}\Bigr|=o(1),\qquad |z|\to 1,
$$
and hence, for small $\ve$, for fixed $R''$, and for $|z|$ close to $1$
we have
\begin{multline}
\card \bigl(\Gamma\cap \mathcal D(w,(R+\ve^3)\rho(z))\bigr)\ge
\card \Gamma(w,R)\\
\ge (q_{-}(R)-\ve^3)R^2,
\qquad  w\in \mathcal D(z,R''\rho(z)).
\label{n10}
\end{multline}
In the same way,
\begin{multline}
E\subset E'=\bigl\{w\in \mathcal D(z,R''\rho(z)/2): \\
\card \bigl(\Gamma\cap \mathcal
D(w,(R+\ve^3)\rho(z))\bigr) \ge (q_{-}(R)+\delta)R^2\bigr\}.
\label{n11}
\end{multline}

We use that by the Fubini theorem, for $0<r_1<r_2$ and for $F\subset \mathcal D(r_2-r_1)$,
\begin{equation}
\card F=\frac{1}{\pi r_1^2}\int_{\mathcal D(r_2)}\card(F\cap \mathcal D(w,r_1))\, dm_2(w).
\label{n12}
\end{equation}

(i) By \eqref{n12},
\begin{multline}
\frac{\card \Gamma(z,R'')}{R^{\prime\prime2}}\\
\ge\int_{\mathcal D(z, (R''-R-\ve^3)\rho(z))}
\!\!\!\!
\frac{
\card\bigl(\Gamma\cap \mathcal
D(w,(R+\ve^3)\rho(z))\bigr)}
{\pi(R+\ve^3)^2\rho(z)^2R^{\prime\prime2}}
\,dm_2(w).\label{sp3}
\end{multline}
Therefore, by \eqref{n10}, for $|z|$ close to $1$,
\begin{multline*}
\frac{\card \Gamma(z,R'')}{R^{\prime\prime2}}
\\ \ge
(q_{-}(R)-\ve^3)\Bigl(\frac{R^{\prime\prime}-R-\ve^3}{R^{\prime\prime}}\Bigr)^2
\Bigl(\frac{R}{R+\ve^3}\Bigr)^2
\ge q_{-}(R)-\ve,
\end{multline*}
for $\ve<\ve(R)$, $R''\ge R'(R,\varepsilon)$.

(ii) By \eqref{n10}--\eqref{n11}, for small $\ve>0$, fixed $R''$ and $|z|$ close to $1$ we have
\begin{gather*}
\int_{\mathcal D(z, (R''-R-\ve^3)\rho(z))}
\frac
{\card\bigl(\Gamma\cap \mathcal D(w,(R+\ve^3)\rho(z))\bigr)}
{\pi(R+\ve^3)^2\rho(z)^2R^{\prime\prime2}} dm_2(w)\\
=\int_{E'}\ldots + \int_{\mathcal D(z, (R''-R-\ve^3)\rho(z))\setminus E'}
\\ 
\ge
\frac{R^2m_2E'}{\pi(R+\ve^3)^2\rho(z)^2R^{\prime\prime2}}
(q_{-}(R)+\delta)\\ +
\frac{R^2(\pi(R''-R-\ve^3)^2\rho(z)^2- m_2E')}{\pi(R+\ve^3)^2\rho(z)^2R^{\prime\prime2}}
\bigl(q_{-}(R)-\ve^3\bigr).
\end{gather*}
If $m_2E'\ge \delta m_2\mathcal D(z,R''\rho(z)/2)$, then by
\eqref{sp3} we obtain
\begin{multline*}
\frac{\card \Gamma(z,R'')}{R^{\prime\prime2}}\\
\ge\int_{\mathcal D(z, (R''-R-\ve^3)\rho(z))}
\frac
{\card\bigl(\Gamma\cap \mathcal D(w,(R+\ve^3)\rho(z))\bigr)}
{\pi(R+\ve^3)^2\rho(z)^2R^{\prime\prime2}} dm_2(w)
\\
\ge \Bigl(\frac{R}{R+\ve^3}\Bigr)^2\Bigl[\Bigl(\frac{R''-R-\ve^3}{R''}\Bigr)^2
(q_{-}(R)-\ve^3)+\frac{\delta}4(\delta+\ve^3)\Bigr]\\
\ge q_{-}(R)+\frac{\delta^2}{5}
\end{multline*}
for $\ve=\ve(\delta)$, $R''\ge R'(R,\delta)$, $|z|\ge\eta_2(R'')$.
\end{proof}

By Lemma~\ref{subl}(i), for every $R_0$ and $\ve$ such that $0<\ve<\ve(R_0)$, we have
$$
D^{-}_\rho(\Gamma,\D)=\liminf_{R\to\infty} q_{-}(R)\ge q_{-}(R_0)-\ve.
$$
Therefore, we obtain

\begin{cor}
\begin{equation}
\lim_{R\to\infty} q_{-}(R)=D^{-}_\rho(\Gamma,\D),
\label{n21}
\end{equation}
and
\begin{equation}
q_{-}(R)\le D^{-}_\rho(\Gamma,\D), \qquad R> 0.
\label{n20}
\end{equation}
\label{m7}
\end{cor}

Given closed subsets $A$ and $B$ of $\CC$, the Fr\'echet distance
$[A,B]$ is the smallest $t>0$ such that $A\subset
B+t\overline{\D}$, $B\subset A+t\overline{\D}$. A sequence
$\{A_n\}$, $A_n\subset\CC$, {\it converges weakly} to $A\subset\CC$ if
for every $R>0$,
$$
\bigl[(A_n\cap \mathcal D(R))\cup R\T,(A\cap \mathcal D(R))\cup
R\T\bigr]\to 0, \qquad n\to\infty.
$$
In this case we use the notation $A_n\rightharpoonup A$. Given any
sequence $\{A_n\}$, $A_n\subset\CC$, we can choose a weakly
convergent subsequence $\{A_{n_k}\}$.

\begin{lem} If $\Gamma\subset\D$, and
$D^{-}_\rho(\Gamma,\D)\le \frac 12$, then there exists a sequence
of points $z_j\in\D$, $|z_j|\to 1$, a sequence $R_j\to\infty$, $j\to\infty$, and a  
subset $\Gamma_0$ of $\CC$ such that
\begin{gather}
\Gamma^{\#}(z_j,R_j)\rightharpoonup\Gamma_0, \qquad j\to\infty, \label{17}\\
\liminf_{R\to\infty}\frac{\card(\Gamma_0\cap \mathcal
D(R))}{R^2}\le \frac 12, \label{18}
\end{gather}
where
\begin{multline*}
\Gamma^{\#}(z,R)=\bigl\{w\in\CC:z+w\rho(z)\in\Gamma(z,R)\bigr\}\\
=\bigl\{w\in\mathcal D(R):z+w\rho(z)\in\Gamma\bigr\}.
\end{multline*}
\label{ll9}
\end{lem}

\begin{proof} Choose a sequence of positive numbers 
$\delta_k$,
$\sum_{k\ge 1}\delta_k\le 1$, set $r_k=2^k$, $k\ge 1$,
and apply Lemma~\ref{subl}(ii) to find $\ve_1>0$,
$0<\eta_1<1$, $R_1$ such that for $\eta_1\le|w|<1$, if
$$
\frac{\card \Gamma(w,R_1)}{R_1^2}\le q_{-}(R_1)+\ve_1,
$$
then there exists $z=z_1(w)\in\mathcal D(w,R_1\rho(w)/2)$
such that
$$
\frac{\card \Gamma(z,r_1)}{r_1^2}\le q_{-}(r_1)+\delta_1.
$$ 

Applying Lemma~\ref{subl}(ii) repeatedly, we find 
$\ve_m\to 0$, $\eta_m\to 1$, $R_m\to\infty$, $m\to\infty$,
such that for $m\ge 1$, $\eta_m\le|w|<1$, if
$$
\frac{\card \Gamma(w,R_m)}{R_m^2}\le q_{-}(R_m)+\ve_m,
$$
then there exists $z=z_m(w)\in\mathcal D(w,R_m\rho(w)/2)$,
such that
$$
\frac{\card \Gamma(z,r_k)}{r_k^2}\le q_{-}(r_k)+\delta_k,\qquad 1\le k\le m.
$$
Next, by the definition of $q_-(R_m)$, we can find $w_m\in\D$,
$\eta_m\le |w_m|<1$, such that
$$
\frac{\card \Gamma(w_m,R_m)}{R_m^2}\le q_-(R_m)+\ve_m,
$$
and define $z_m=z_m(w_m)$. We obtain
\begin{equation}
\limsup_{m\to\infty}\frac{\card(\Gamma(z_m,r_k))}{r_k^2}\le
D^{-}_\rho(\Gamma,\D)\le \frac12,\qquad k\ge 1. \label{19}
\end{equation}
Finally, we choose a sequence $\{m_k\}$ and a set $\Gamma_0\in\CC$
such that
$$
\Gamma^{\#}(z_{m_k},r_{m_k})\rightharpoonup\Gamma_0,\qquad
m_k\to\infty.
$$
The property \eqref{18} follows from \eqref{19}.
\end{proof}

Analogously, we have

\begin{lem} If $\Gamma\subset\D$, and
$D^{+}_\rho(\Gamma,\D)\ge \frac 12$, then there exists a sequence
of points $z_j\in\D$, $|z_j|\to 1$, a sequence $R_j\to\infty$, $j\to\infty$, and a  
subset $\Gamma_0$ of $\CC$ such that
\begin{gather}
\Gamma^{\#}(z_j,R_j)\rightharpoonup\Gamma_0, \qquad j\to\infty, \label{17b}\\
\limsup_{R\to\infty}\frac{\card(\Gamma_0\cap \mathcal D(R))}{R^2}\ge \frac 12. \label{18b}
\end{gather}
\label{ll9b}
\end{lem}

\section{Sampling theorems}
\label{ssa}

We set $\beta(z)=|z|^2/4$.

\begin{proof}[Proof of Theorem~{\rm \ref{tt01}}]
By Corollary~\ref{m2}, every sampling set for $\ahd$
contains a $d_\rho$-separated subset which is also a sampling set for $\ahd$.

(A) Suppose that $\Gamma$ is $d_\rho$-separated, and
$D^{-}_\rho(\Gamma,\D)\le \frac 12$. We follow the scheme proposed
in \cite{LS}. We apply Lemma~\ref{ll9} to obtain $z_j$, $R_j$, and
$\Gamma_0$ satisfying \eqref{17}--\eqref{18}. Fix $\varepsilon>0$.
By the theorem of Seip on sampling in Fock type spaces
\cite[Theorem~2.3]{S}, there exists $f\in
{\mathcal{A}_{\beta}(\CC)}$ such that
$$
\|f\|_{\beta}=1,\qquad \|f\|_{\beta,\Gamma_0}\le \varepsilon.
$$
For $K>1$ we set $f_K(z)=f((1-K^{-3/2})z)$. Then
\begin{multline*}
|f_K(z)|e^{-\beta(z)}\le |f_K(z)|e^{-(1-K^{-3/2})^2\beta(z)}\\
\le
|f(z)|e^{-\beta(z)}+\bigl||f(z)|e^{-\beta(z)}-|f((1-K^{-3/2})z)|
e^{-(1-K^{-3/2})^2\beta(z)}\bigr|
\\=|f(z)|e^{-\beta(z)}+o(1),\qquad |z|\le K,\quad K\to\infty,
\end{multline*}
where in the last relation we use \cite[Lemma~3.1]{S} (for a similar estimate see Lemma~\ref{m1}~(i)).
Furthermore,
$$
|f_K(z)|e^{-\beta(z)}=o(1),\qquad |z|> K,\quad K\to\infty.
$$
Therefore, for sufficiently large $K$ we get
$$
\|f_K\|_{\beta}\asymp 1,\qquad \|f_K\|_{\beta,\Gamma_0}\le 2\varepsilon.
$$
We fix such $K$ and for $N\ge 0$ set
$$
T_Nf_K(z)=\sum_{0\le n\le N}c_nz^n,
$$
where
$$
f_K(z)=\sum_{n\ge 0}c_nz^n.
$$
As in \cite[page 169]{LS}, by the Cauchy formula,
$$
|c_n|\le c\cdot \inf_r\frac{\exp[(1-K^{-3/2})^2r^2/4]}{r^n},
$$
and hence,
\begin{equation}
\sum_{n\ge 0}|c_n z^n|e^{-\beta(z)}\le c(1+|z|)^4e^{[(1-K^{-3/2})^2-1]|z|^2/4},
\qquad z\in\CC.
\label{fd143}
\end{equation}
Therefore, for sufficiently large $N$ we have
$$
\|T_Nf_K\|_{\beta}\asymp 1,\qquad \|T_Nf_K\|_{\beta,\Gamma_0}\le
3\varepsilon.
$$
We fix such $N$, set $P=T_Nf_K$, and choose $a\in\CC$ such that
$$
|P(a)|e^{-\beta(a)}\asymp 1.
$$

By \eqref{17}, we can find large $R>|a|$ and $z$ close to the unit circle such that
\begin{gather*}
|P(w)|\le\ve |w|^R, \qquad |w|\ge R,\\
\|P\|_{\beta,\Gamma^{\#}(z,R)}\le 4\varepsilon.
\end{gather*}
We set $z^*=z+a\rho(z)$,
apply Proposition~\ref{ll7} to get $g=g_{z,R}$, and define
$$
f(w)=g(w)P\Bigl(\frac{w-z}{\rho(z)}\Bigr).
$$
Then $f\in\ahd$,
\begin{gather*}
|f(z^*)|e^{-h(z^*)}\asymp |g(z^*)|e^{-h(z^*)}\cdot|P(a)|\asymp
e^{-|a|^2/4}e^{\beta(a)}=1,\\
|f(w)|e^{-h(w)}\le e^{-|w-z|^2/(4\rho(z)^2)}\cdot 4\ve
e^{|w-z|^2/(4\rho(z)^2)} =4\ve,\quad w\in\Gamma(z,R),
\end{gather*}
\begin{multline*}
|f(w)|e^{-h(w)}\le
c\Bigl[\frac{R^2\rho(z)^2}{e|z-w|^2}\Bigr]^{R^2/4}\cdot
\ve \Bigl[\frac{|z-w|}{\rho(z)}\Bigr]^{R}\\
\le c\ve,\qquad w\in\mathbb D\setminus \mathcal D(z,R\rho(z)),
\end{multline*}
with $c$ independent of $\varepsilon,R$. Since $\varepsilon$ can be chosen arbitrarily small,
this shows that $\Gamma$ is not a sampling set for $\ahd$.

(B) Now we assume that $\Gamma$ is a $d_\rho$-separated subset of $\D$, $D^{-}_\rho(\Gamma,\D)> \frac 12$,
and $\Gamma$ is not a sampling set for $\ahd$. Then there exist functions $f_n\in\ahd$
such that $\|f_n\|_{h}=1$ and $\|f_n\|_{h,\Gamma}\to 0$, $n\to\infty$.
By the normal function argument, either (B1) $f_n$ tend to $0$
uniformly on compact subsets of the unit disc or (B2)
there exists a subsequence $f_{n_k}$ converging uniformly on compact subsets of
the unit disc to $f\in\ahd$,
$f\ne 0$, with $f\bigm|\Gamma =0$.

In case (B1), using Proposition~\ref{ll7}
we can find $z_n\in\D$, $R_n\to \infty$, $n\to\infty$,
such that the functions
$$
F_n(w)=\frac{f(z_n+w\rho(z_n))}{g_{z_n,R_n}(z_n+w\rho(z_n))}
$$
satisfy the conditions:
\begin{align*}
|F_n(0)|&\asymp 1,\\
|F_n(w)|&\le e^{|w|^2/4},\qquad |w|\le R_n,\\
\sup_{\Gamma^{\#}(z_n,R_n)}|F_n|&\to 0,\qquad n\to\infty.
\end{align*}
By Corollary~\ref{m7}, we can find $q$, $\frac 12<q<\dmGD$,
and $0<C<R'_n<R_n$, $R'_n\to\infty$ as $n\to\infty$, such that
$$
\card \Gamma^{\#}(z_n,r)\ge qr^2,\qquad C\le r\le R'_n.
$$
Again by the normal function argument, we can choose a sequence
$n_k\to\infty$, $k\to\infty$, such that $F_{n_k}$ converge
uniformly on compact subsets of $\CC$ to
$F\in{\mathcal{A}_{\beta}(\CC)}$, and
$\Gamma(z_{n_k},R_{n_k})\rightharpoonup\Gamma^*$ such that
\begin{align*}
F(0)&\ne 0,\\
F\bigm |\Gamma^*&=0,\\
\card (\Gamma^*\cap \clos\mathcal D(r))&\ge qr^2,\qquad r\ge C.
\end{align*}
To get the last inequality we use that $\Gamma$ is 
$d_\rho$-separated.

However, by Jensen's inequality,
\begin{gather*}
-\infty< \log |F(0)|\le \frac{1}{2\pi}\int^{2\pi}_0 \log |F(re^{i\theta})|d\theta
-\sum_{w_k\in\Gamma^*\cap \mathcal D(r)}\log \frac{r}{|w_k|}\\
\le \frac{r^2}4-
\int^r_0\log \frac rs\,dn(s)=\frac{r^2}4-
\int^r_0\frac {n(s)}{s}\,ds\\
\le \frac{r^2}4- \frac{qr^2}2+O(1)\to -\infty,\qquad r\to\infty,
\end{gather*}
where $n(r)=\card (\Gamma^*\cap \clos\mathcal D(r))$.
This contradiction implies our assertion in case (B1).

In case (B2),
without loss of generality we can assume that $0\not\in\Gamma$, $f(0)\ne 0$. By Jensen's inequality,
\begin{equation}
-\infty< \log |f(0)|\le \frac{1}{2\pi}\int^{2\pi}_0 \log |f(re^{i\theta})|d\theta
-\sum_{w_k\in\Gamma\cap \mathcal D(r)}\log \frac{r}{|w_k|}.
\label{n26}
\end{equation}
Furthermore, we choose $\ve>0$ and large $R>R(\ve)$.
Then
\begin{multline*}
\sum_{w_k\in\Gamma\cap \mathcal D(r)}\log \frac{r}{|w_k|}\\ =
\sum_{w_k\in\Gamma\cap \mathcal D(r)} \log \frac{r}{|w_k|}\cdot
\frac{1}{\pi R^2\rho(w_k)^2}\int_{\mathcal D(w_k,R\rho(w_k))}dm_2(w)\\ \ge
O(1)+(1-\ve)\int_{\mathcal D(r-R^2\rho(r))}\frac{1}{\pi R^2\rho(w)^2}\cdot
\log \frac{r}{|w|}\times\\
\times\card\bigl[w_k\in \Gamma : w\in\mathcal
D(w_k,R\rho(w_k))\bigr]\,dm_2(w)\\ \ge
O(1)+\frac{(1-\ve)^2q_-(R-\ve)}\pi \int_{\mathcal
D(r-R^2\rho(r))}\frac{1}{\rho(w)^2}\log \frac{r}{|w|} \,dm_2(w)\\
=O(1)+\frac{(1-\ve)^2q_-(R-\ve)}\pi 2\pi h(r),\qquad r\to 1,
\end{multline*}
that contradicts to \eqref{n26} for small $\ve>0$ and $R>R(\ve)$. This proves our assertion.
\end{proof}

\begin{proof}[Proof of Theorem~{\rm \ref{tt01b}}]
By Lemmas~\ref{m4} and ~\ref{m6}, every sampling set for $\ahdp$
is a finite union of $d_\rho$-separated subsets and
contains a $d_\rho$-separated subset which is also a sampling set for $\ahdp$.

(A) Suppose that $\Gamma$ is a $d_\rho$-separated subset of the
unit disc and $D^{-}_\rho(\Gamma,\D)\le \frac 12$. As in part (A)
of the proof of Theorem~\ref{tt01}, we apply Lemma~\ref{ll9} to
obtain $z_j$, $R_j$, and $\Gamma_0$ satisfying
\eqref{17}--\eqref{18}. Furthermore, $\Gamma_0$ is uniformly
separated, that is
$$
\inf\{|z_1-z_2|:z_1,z_2\in\Gamma_0,\, z_1\ne z_2\}>0.
$$
Then, by a version of a result of Seip \cite[Lemma~7.1]{S} (see also Lem\-ma~\ref{m4}),
\begin{equation}
\sum_{z\in\Gamma_0}e^{-p\beta(z)}|g(z)|^p\le c\|g\|^p_{p,\beta},
\qquad g\in\mathcal{A}^p_\beta(\CC).
\label{n19}
\end{equation}

Fix $\varepsilon>0$. By the theorem of Seip on sampling in Fock
spaces \cite[Theorem~2.1]{S}, there exists $f\in
{\mathcal{A}^p_{\beta}(\CC)}$ such that
$$
\|f\|_{p,\beta}=1,\qquad \|f\|_{p,\beta,\Gamma_0}\le \varepsilon.
$$
We approximate $f$ by a polynomial $P$ in the norm of
$\mathcal{A}^p_\beta(\CC)$ and obtain, using \eqref{n19}, that
$$
1-\varepsilon\le \|P\|_{p,\beta}\le 1+\varepsilon ,\qquad
\|P\|_{p,\beta,\Gamma_0}\le 2\varepsilon.
$$
For some $M>0$ we have
\begin{equation}
\int_{\mathcal D(M)}|P(z)|^pe^{-p\beta(z)}dm_2(z)\ge\frac12.
\label{n37f}
\end{equation}
By \eqref{17}, we can find large $R>M$ and $z$
 close to the unit circle such that
\begin{gather}
|P(w)|\le\ve |w|^R, \qquad |w|\ge R.\label{sp4}\\
\|P\|_{p,\beta,\Gamma^{\#}(z,R)}\le 3\varepsilon,\label{n34}
\end{gather}

We apply Proposition~\ref{ll7} to get $g=g_{z,R}$, and define
$$
f(w)=\frac{g(w)}{\rho(z)^{2/p}}P\Bigl(\frac{w-z}{\rho(z)}\Bigr).
$$
Then, by \eqref{15} and \eqref{n37f},
\begin{multline}
\|f\|^p_{p,h}\ge \int_{|w-z|\le R\rho(z)}|f(w)|^pe^{-ph(w)}dm_2(w)\\ \asymp
\int_{|w-z|\le R\rho(z)}\frac 1{\rho(z)^2}\Bigl|P\Bigl(\frac{w-z}{\rho(z)}\Bigr)\Bigr|^p
e^{-p|w-z|^2/[4\rho(z)^2]}dm_2(w)
\\ =\int_{\mathcal D(R)}|P(w)|^pe^{-p|w|^2/4}dm_2(w)\ge \frac 12.
\label{n36}
\end{multline}

On the other hand,
\begin{equation*}
\|f\|^p_{p,h,\Gamma}=\sum_{w\in\Gamma}|f(w)|^pe^{-ph(w)}\rho(w)^2=
\sum_{w\in\Gamma(z,R)}\ldots+\sum_{w\in\Gamma\setminus\Gamma(z,R)}\ldots\,.
\end{equation*}
By \eqref{15} and \eqref{n34},
\begin{multline*}
\sum_{w\in\Gamma(z,R)}\frac{\rho(w)^2}{\rho(z)^2}|g(w)|^pe^{-ph(w)}
\Bigl|P\Bigl(\frac{w-z}{\rho(z)}\Bigr)\Bigr|^p\\
\le
c\sum_{w\in\Gamma^{\#}(z,R)}e^{-p|w|^2/4}|P(w)|^p\le C\ve^p.
\end{multline*}
Since $\Gamma$ is $d_\rho$-separated, by Lemma~\ref{m3} we have
\begin{equation*}
\sum_{w\in\Gamma\setminus\Gamma(z,R)}|f(w)|^pe^{-ph(w)}\rho(w)^2\le
c\int_{\D\setminus \mathcal D(z,(R-1)\rho(z))}
\!\!|f(w)|^pe^{-ph(w)}dm_2(w).
\end{equation*}
Furthermore, by \eqref{15}, \eqref{16}, and \eqref{sp4},
\begin{gather*}
\int_{\D\setminus \mathcal D(z,(R-1)\rho(z))}
\frac 1{\rho(z)^2}|g(w)|^pe^{-ph(w)}\Bigl|P\Bigl(\frac{w-z}{\rho(z)}\Bigr)\Bigr|^pdm_2(w)\\
\le c\ve^p\int_{\D\setminus \mathcal D(z,R\rho(z))}
\frac 1{\rho(z)^2}\Bigl[\frac{|w-z|}{\rho(z)}\Bigr]^{pR}
\Bigl[\frac{R^2\rho(z)^2}{e|w-z|^2}\Bigr]^{pR^2/4} dm_2(w)\\
+c\ve^p\int_{\mathcal D(z,R\rho(z))\setminus \mathcal D(z,(R-1)\rho(z))  }
\frac 1{\rho(z)^2}R^{pR}
e^{-p|w-z|^2/[4\rho(z)^2]}dm_2(w)\\ =
c\ve^p\int_{|w|>R}|w|^{pR}\Bigl[\frac{R^2}{e|w|^2}\Bigr]^{pR^2/4} dm_2(w)\\
+c\ve^p\int_{R-1<|w|<R}R^{pR}e^{-p|w|^2/4} dm_2(w)
\le c \ve^p,
\end{gather*}
with $c$ independent of $R$.
This together with \eqref{n36} shows that $\Gamma$ is not a sampling set for $\ahdp$.

(B) Now we assume that $\Gamma$ is a $d_\rho$-separated subset of $\D$,
$D^{-}_\rho(\Gamma,\D)> \frac 12$.
Then, by Theorem~\ref{tt01}, we can fix small $\ve>0$
such that $\Gamma$ is a sampling set
for $\mathcal{A}_{(1+\ve)h}(\D)$.
Following the method of \cite[Section 6]{BS}, we 
are going to prove that $\Gamma$ is a sampling set
for $\ahdp$.

We set
\begin{gather*}
\mathcal{A}_{(1+\ve)h,0}(\D)=\bigl\{F\in\hol(\D):\lim_{|z|\to 1}F(z)e^{-(1+\ve)h(z)}=0\bigr\},\\
R_\Gamma: F\in \mathcal{A}_{(1+\ve)h,0}(\D) \mapsto
\bigl\{F(z_k)e^{-(1+\ve)h(z_k)}\bigr\}_{z_k\in\Gamma}\in c_0.
\end{gather*}
Since $\Gamma$ is a sampling set for $\mathcal{A}_{(1+\ve)h}(\D)$,
$R_\Gamma$ is an invertible linear operator onto a closed subspace $V$
of the space $c_0$. Therefore, linear functionals $E_z: v\in V\mapsto (R_\Gamma^{-1}v)(z)e^{-(1+\ve)h(z)}$,
$\|E_z\|\asymp1$, $z\in\D$, (see also Lemma~\ref{m8}) extend to linear functionals on
$c_0$ bounded uniformly in $z$. Thus, for every $z\in\D$, there exist $b_k(z)$, $k\ge 1$,
such that
\begin{equation}
\sum_{k\ge 1} |b_k(z)|\le C,
\label{52}
\end{equation}
with $C$ independent of $z$, and
\begin{equation}
F(z)e^{-(1+\ve)h(z)}=\sum_{k\ge 1} b_k(z)F(z_k)e^{-(1+\ve)h(z_k)}.
\label{n28}
\end{equation}

Let $f\in\ahdp$. By Lemma~\ref{m9}, $f\in
\mathcal{A}_{(1+\ve)h,0}(\D)$. For every $z\in\D$, we use
Lemma~\ref{m10} (with $h$ replaced by $\ve h$) to get $g_z$
satisfying \eqref{n30}--\eqref{n31} (with $\rho$ replaced by $\varepsilon^{-1/2}\rho$). 
Now, we apply \eqref{n28} to
$F=fg_z$. By \eqref{n30} we obtain
$$
|f(z)|e^{-h(z)}\le C\sum_{k\ge 1} |b_k(z)||f(z_k)|e^{-h(z_k)}|g_z(z_k)|e^{-\ve h(z_k)}.
$$
By \eqref{52},
$$
|f(z)|^pe^{-ph(z)}\le C\sum_{k\ge 1} |f(z_k)|^pe^{-ph(z_k)}|g_z(z_k)|^pe^{-p\ve h(z_k)}
$$
and hence,
$$
\|f\|^p_{p,h} \le C\sum_{k\ge 1} 
|f(z_k)|^pe^{-ph(z_k)}\int_\D |g_z(z_k)|^pe^{-p\ve h(z_k)}dm_2(z).
$$
It remains to note that by \eqref{n31},
\begin{multline*}
\int_\D |g_z(z_k)|^pe^{-p\ve h(z_k)}dm_2(z)\\ \le c(\ve) \int_\D\min\Bigl[1,\frac{\rho(z_k)^{3p}}{|z-z_k|^{3p}}\Bigr]\,dm_2(z) \le
C(\ve)\rho(z_k)^2. 
\end{multline*}
Then $\|f\|_{p,h} \le C\|f\|_{p,h,\Gamma}$, and using
Lemma~\ref{m4}, we conclude that $\Gamma$ is a sampling set for $\ahdp$.
\end{proof}
\smallskip

\section{Interpolation theorems}
\label{sin}

\begin{proof}[Proof of Theorem~{\rm \ref{tti01}}]
Corollary~\ref{m15} claims that every interpolation set for $\ahd$
is $d_\rho$-separated.

(A) Let $\Gamma$ be a $d_\rho$-separated subset of $\D$, and
let $D^{+}_\rho(\Gamma,\D)\ge \frac 12$.
Suppose that $\Gamma$ is an interpolation set for $\ahd$.
We apply Lemma~\ref{ll9b} to obtain $z_j$, $R_j$, and $\Gamma_0=\{z^0_k\}_{k\ge 1}$
satisfying \eqref{17b}, \eqref{18b}. Suppose that $a_k\in \CC$ satisfy the estimate
$$
|a_k|e^{-\beta(z^0_k)}\le 1.
$$
By \eqref{17b}, we can choose a subsequence $\{z'_j\}$ of $\{z_j\}$ and $R'_j\to\infty$
satisfying the following properties: $\Gamma(z'_j,R'_j)$ are disjoint,
$$
B_j=\card \Gamma(z'_j,R'_j)= \card [\Gamma_0\cap \mathcal
D(R'_j)],
$$
and we can enumerate $\Gamma(z'_j,R'_j)=\{w_{jk}\}_{1\le k\le B_j}$ in such a way that
\begin{equation}
\max_{1\le k\le B_j}
\Bigl|z^0_k-\frac{w_{jk}-z'_j}{\rho(z'_j)}\Bigr|\to 0,\qquad j\to \infty.
\label{n41}
\end{equation}
Without loss of generality, we can assume that $|z'_j|>1-\eta(R'_j)$. Therefore,
by Proposition~\ref{ll7}, there exist $g_j=g_{z'_j,R'_j}$ satisfying \eqref{15}
and \eqref{16}.

For $j\ge 1$ we consider the following interpolation problem:
\begin{equation}
f_j(w)=
\begin{cases}
a_kg_j(w_{jk}), & w=w_{jk}\in\Gamma(z'_j,R'_j),\\
0, & w\in\Gamma\setminus\Gamma(z'_j,R'_j).
\end{cases}
\label{n42}
\end{equation}
By our assumption on $\Gamma$, we can find $f_j\in\ahd$ satisfying \eqref{n42}.
Then the functions
$$
F_j(w)=\frac{f_j(z'_j+w\rho(z'_j))}{g_j(z'_j+w\rho(z'_j))}
$$
satisfy the properties
\begin{gather*}
F_j\Bigl(\frac{w_{jk}-z'_j}{\rho(z'_j)}\Bigr)=a_k,\qquad 1\le k\le B_j,\\
|F_j(w)|\le ce^{\beta(w)},\qquad |w|\le R'_j.
\end{gather*}

By a normal families argument and by \eqref{n41}, we conclude that there exists
an entire function $F\in \mathcal A_\beta(\CC)$ such that
$$
F(z^0_k)=a_k,\qquad k\ge 1.
$$
Thus, $\Gamma_0$ is a set of interpolation for $\mathcal A_\beta(\CC)$.
However, by the theorem of Seip on interpolation in the Fock type spaces
\cite[Theorem~2.4]{S}, this is impossible for $\Gamma_0$ satisfying \eqref{18b}.
This contradiction proves our assertion.

(B) Now we assume that $\Gamma$ is a $d_\rho$-separated subset of $\D$,
$D^{+}_\rho(\Gamma,\D)< \frac 12$.
First of all, if $\Gamma$ is an interpolation set for $\ahd$, and
$\lambda\in\D\setminus\Gamma$, then $\Gamma\cup\{\lambda\}$ is also
an interpolation set for $\ahd$. (Later on, to deal with the plane case, we add to $\Gamma$ an infinite sequence in such a way that the modified $\Gamma$ satisfies the same conditions,
and then use that if $\Gamma$ is an interpolation set for $\ahc$, and
$\lambda\in\CC\setminus\Gamma$, $\mu\in\Gamma$, then $\Gamma\cup\{\lambda\}\setminus\{\mu\}$ is also
an interpolation set for $\ahc$.)

For every sufficiently large $R$, $R\ge R_0$, we can find $\eta_1(R)<1$ such that the family
of sets $\Gamma^{\#}(z,R)$, $z\in\Gamma\setminus \mathcal D(\eta_1(R))$, $R\ge R_0$,
satisfies the uniform estimates
\begin{equation}
\left.
\begin{gathered}
\sup_{w\in\CC,\, r\ge R_0}\frac{\card\bigl(\mathcal D(w,r)\cap \Gamma^{\#}(z,R)\bigr)}{r^2}
<\frac 12,\\
\inf_{w_1\ne w_2,\,w_1,w_2\in \Gamma^{\#}(z,R)}|w_1-w_2|>0.
\end{gathered}
\right\}\label{f86}
\end{equation}
Therefore, by (a variant of) the theorem of Seip-Wallst\'en on interpolation in the Fock
type spaces \cite[Theorem~2.4]{S}, \cite[Theorem~1.2]{SW}, for some $c<\infty$, $\ve>0$,
and for every $R\ge R_0$, $z\in\Gamma\setminus \mathcal
D(\eta_1(R))$, there exists $F_{z,R}\in A_{(1-\ve)\beta}(\CC)$,
such that
\begin{equation}
\left.
\begin{gathered}
F_{z,R}(0)=1,\\
F_{z,R}\bigm| \Gamma^{\#}(z,R)\setminus\{0\}=0,\\
\|F_{z,R}\|_{(1-\ve)\beta}\le c.
\end{gathered}
\right\}
\label{n38}
\end{equation}

To continue, we need a simple estimate similar to \eqref{fd143}.

\begin{lem} If $F(z)=\sum_{n\ge 0}c_nz^n$,
$|F(z)|\le \exp|z|^2$, $z\in\CC$, $N\in\Z_+$, and
$$
(T_NF)(z)=\sum_{0\le n\le N}c_nz^n,
$$
then
\begin{gather*}
|(F-T_NF)(z)|\le 2^{-(N-3)/2},\qquad |z|\le \sqrt{N/(4e)},\\
|(T_NF)(z)|\le (N+1)\exp |z|^2,\qquad \sqrt{N/(4e)}< |z|\le \sqrt{N/2},\\
|(T_NF)(z)|\le (N+1)|z|^N(2e/N)^{N/2},\qquad |z|>\sqrt{N/2}.
\end{gather*}
\label{m22}
\end{lem}

\begin{proof} By the Cauchy formula, we have
$$
|c_n|\le \inf_{r>0}\bigl[r^{-n}\exp r^2\bigr]=\exp\Bigl[-\frac
n2\log\frac n{2e}\Bigr],
$$
and hence,
$$
\sum_{n>N}|c_n|r^n\le
\sum_{n>N}\Bigl(\frac{n}{2e}\cdot\frac{4e}{N}\Bigr)^{-n/2}\le
2^{-(N-3)/2}, \qquad r\le \sqrt{N/(4e)},
$$
and
$$
\sum_{0\le n\le N}|c_n|r^n\le (N+1)\max_{0\le n\le N}
\exp\Bigl[-\frac n2\log\frac n{2e}+n\log r\Bigr],\qquad r\ge 0.
$$
Furthermore,
$$
-\frac n2\log\frac n{2e}+n\log r\le r^2, \qquad r\ge 0,
$$
and
$$
-\frac n2\log\frac n{2e}+n\log r\le N\log r -\frac N2\log \frac {N}{2e},
\quad r\ge\sqrt{N/2},\,0\le n\le N.
$$
\end{proof}

\begin{cor} If $\ve>0$, $F\in \mathcal A_{(1-\ve)\beta}(\CC)$, $\|F\|_{(1-\ve)\beta}\le 1$,
$N\in\Z_+$, 
$T_NF$ is defined as above, and $R=\sqrt{2N/(1-\ve)}$ is
sufficiently large, $R>R(\ve)$, then for some $c=c(\varepsilon)>0$
independent of $z,R$ we have
\begin{align*}
|(F-T_NF)(z)|e^{-|z|^2/4}&\le e^{-cR^2},\qquad |z|\le R,\\
|(T_NF)(z)|e^{-|z|^2/4}&\le 2e^{-c|z|^2},\qquad |z|\le R,\\
|(T_NF)(z)|\Bigl(\frac{e|z|^2}{R^2}\Bigr)^{-R^2/4}&\le
\Bigl(\frac{e|z|^2}{R^2}\Bigr)^{-\ve R^2/5},\qquad |z|> R.
\end{align*}
\label{m17}
\end{cor}

For large $N$ we set $R=\sqrt{2N/(1-\ve)}$, and for $z\in\Gamma$
sufficiently close to the unit circle, define, using $g_{z,R}$
from Proposition~\ref{m14} and $F_{z,R}$ from \eqref{n38},
\begin{equation}
U_z(w)=g_{z,R}(w)\cdot (T_NF_{z,R})\Bigl(\frac{w-z}{\rho(z)}\Bigr).
\label{47a}
\end{equation}
If $\Gamma=\{z_n\}_{n\ge 1}$,
$a_n\in\CC$, $n\ge 1$, and
\begin{equation}
\sup_{n\ge 1} |a_n|e^{-h(z_n)}\le 1,
\label{f31}
\end{equation}
then for $|z_n|\ge \eta=\max(\eta(R),\eta_1(R))$ 
($\eta(R)$ is introduced in Proposition~\ref{m14})
we put
\begin{equation}
V_n=\frac{a_nU_{z_n}}{U_{z_n}(z_n)}.
\label{47b}
\end{equation}
Then
\begin{equation}
V_n(z_n)=a_n,\label{40a}
\end{equation}
and by the estimates in Proposition~\ref{m14} and in
Corollary~\ref{m17} we obtain for $|z_n|\ge \eta$ that
\begin{gather}
|V_n(z_k)|e^{-h(z_k)}\le c_0e^{-cR^2},\qquad |z_k-z_n|\le R\rho(z_n),\quad k\ne n,\label{40b}\\
|V_n(z)|e^{-h(z)}\le c_0e^{-c |z-z_n|^2/[\rho(z_n)^2]},\qquad
|z-z_n|\le R\rho(z_n),\label{40c}
\end{gather}
\begin{multline}
|V_n(z)|e^{-h(z)}\le
c_0\Bigl(\frac{R^2\min[\rho(z_n),\rho(z)]^2}{e|z-z_n|^2}\Bigr)^{\ve
R^2/5},\\  |z-z_n|> R\rho(z_n), \label{40d}
\end{multline}
for some $c_0$ independent of $z_n,z,R$.

For $R>1$ we define
$$
A_{R}(z)=\sum_{z_n\in \Gamma,\,
|z-z_n|>R\rho(z_n)}
\Bigl(\frac{R^2\min[\rho(z_n),\rho(z)]^2}{e|z-z_n|^2}\Bigr)^{\ve
R^2/5}.
$$
Suppose that for every $\delta>0$, we can find arbitrarily large $R$ such that
\begin{equation}
\sup_\D A_{R}\le\delta. \label{40e}
\end{equation}
Then for $0<\delta<1/(2c_0)$, for sufficiently large $R$, and for \newline\noindent $\eta=\max(\eta(R),\eta_1(R))$ we can define
$$
f_1=\sum_{z_n\in \Gamma\setminus \mathcal D(\eta)}V_n,
$$
and obtain that for some $B$ independent of $\{a_n\}$ satisfying \eqref{f31},
\begin{gather}
\|f_1\|_{h}\le B,\label{40f}\\
\sup_{z_n\in \Gamma\setminus \mathcal D(\eta)}|f_1(z_n)-a_n|e^{-h(z_n)}\le \frac 12.\label{40g}
\end{gather}
Indeed, by \eqref{40c}, \eqref{40d}, and \eqref{40e}, for $z\in\D$ we have
\begin{multline*}
|f_1(z)|e^{-h(z)}\\
\le \!\!\sum_{z_n\in \Gamma\setminus \mathcal D(\eta),\,|z-z_n|\le R\rho(z_n)}\!\!\!\!|V_n(z)|e^{-h(z)}+
\!\!\sum_{z_n\in \Gamma\setminus \mathcal D(\eta),\, |z-z_n|>R\rho(z_n)}\!\!\!\!|V_n(z)|e^{-h(z)}\\
\le c+c_0\delta.
\end{multline*}
By \eqref{40a}, \eqref{40b}, \eqref{40d}, and \eqref{40e}, for
$z_k\in\Gamma\setminus \mathcal D(\eta)$ we have
\begin{multline*}
|f_1(z_k)-a_k|e^{-h(z)}\\
\le \sum_{|z_k-z_n|\le R\rho(z_n),\, k\ne n}|V_n(z)|e^{-h(z)}+
\sum_{z_n\in \Gamma\setminus \mathcal D(\eta),\, |z_k-z_n|>R\rho(z_n)}|V_n(z)|e^{-h(z)}\\
\le cR^2e^{-cR^2}+c_0\delta\le \frac 12
\end{multline*}
for sufficiently large $R$. We fix such $\delta,R,\eta$.

Iterating the approximation construction and using \eqref{40f} and \eqref{40g},
we obtain $f_2\in\ahd$ such that 
\begin{gather*}
\|f_2\|_{h}\le B/2,\\
\sup_{z_n\in \Gamma\setminus \mathcal D(\eta)}
| f_2(z_n) +f_1(z_n) -a_n|e^{-h(z_n)}\le \frac 14.
\end{gather*}
Continuing this process, we arrive at $f=\sum_{n\ge 1}f_n$ such that
\begin{gather*}
\|f\|_{h}\le 2B,\\
f(z_n)=a_n, \qquad z_n\in \Gamma\setminus \mathcal D(\eta).
\end{gather*}

Thus, $\Gamma\setminus \mathcal D(\eta)$ is a set of interpolation
for $\ahd$, and hence, $\Gamma$ is a set of interpolation for
$\ahd$.

It remains to estimate $A_{R}$ for large $R$. Since $\Gamma$ is $d_\rho$-separated, 
using \eqref{01} we obtain
\begin{gather}
\sum_{z_n\in \Gamma,\,|z-z_n|>R\rho(z_n)}
\Bigl(\frac{R^2\min[\rho(z_n),\rho(z)]^2}{e|z-z_n|^2}\Bigr)^{\ve R^2/5} \notag \\
\le C(\Gamma,h)\int_{\D \setminus \mathcal D(z,R\rho(z))}
\Bigl(\frac{R^2\min[\rho(w),\rho(z)]^2}{e|z-w|^2}\Bigr)^{\ve
R^2/5}\frac{dm_2(w)}{\rho(w)^2}\notag \\ \le
R^2\int_{|\zeta|>1}\Bigl(\frac{1}{e|\zeta|^2}\Bigr)^{\ve
R^2/5}dm_2(\zeta) 
=o(1),\qquad R\to\infty,\label{fd144}
\end{gather}
because for any $z,w\in\D$, $R\ge\sqrt{5/\ve}$,
$$
\Bigl(\frac{\min[\rho(w),\rho(z)]}{\rho(z)}\Bigr)^{2\ve
R^2/5}\Bigl(\frac{\rho(z)}{\rho(w)}\Bigr)^{2}\le 1.
$$
This completes the proof of our assertion.
\end{proof}

\begin{proof}[Proof of Theorem~{\rm \ref{tti01b}}]
By Corollary~\ref{m16}, every set of interpolation for $\ahdp$ is
$d_\rho$-separated.

(A) The argument is analogous to that in the part (A) of
the proof of Theorem~{\rm \ref{tti01}}. We just use
\cite[Theorem~2.2]{S} instead of
\cite[Theorem~2.4]{S}.

(B) Now we assume that $\Gamma$ is a $d_\rho$-separated subset of
$\D$, $D^{+}_\rho(\Gamma,\D)< \frac 12$. As in the part (B) of the
proof of Theorem~{\rm \ref{tti01}} we find $c<\infty$, $\ve>0$,
and $R_0>1$, such that for $R\ge R_0$, $z\in\Gamma\setminus
\mathcal D(\eta_1(R))$, the sets $\Gamma^{\#}(z,R)$
satisfy \eqref{f86}, and there exist $F_{z,R}\in
A^p_{(1-2\ve)\beta}(\CC)$, such that
\begin{equation}
\left.
\begin{gathered}
F_{z,R}(0)=1,\\
F_{z,R}\bigm| \Gamma^{\#}(z,R)\setminus\{0\}=0,\\
\|F_{z,R}\|_{p,(1-2\ve)\beta}\le c.
\end{gathered}
\right\} \label{sp5}
\end{equation}

Instead of Lemma~\ref{m22} and Corollary~\ref{m17} we use

\begin{lem} If $F\in \mathcal A_{p,(1-2\ve)\beta}(\CC)$, 
$\|F\|_{p,(1-2\ve)\beta} \le 1$, if
\newline\noindent
$F(z)=\sum_{n\ge 0}c_nz^n$,
$T_NF$ is defined as in Lemma~{\rm\ref{m22}},
and if $R=\sqrt{2N/(1-\ve)}$ is sufficiently large, $R>R(\ve)$,
then for some $c=c(\ve)>0$ independent of $z,R$ we have
\begin{gather}
|(F-T_NF)(z)|e^{-|z|^2/4}\le e^{-cR^2},\qquad |z|\le R,\label{49b}\\
|(T_NF)(z)|\Bigl(\frac{e|z|^2}{R^2}\Bigr)^{-R^2/4}\le
\Bigl(\frac{e|z|^2}{R^2}\Bigr)^{-\ve R^2/5},\qquad |z|> R,\label{49d}\\
\int_{\mathcal D(R)}|(T_NF)(z)|^pe^{-p|z|^2/4}dm_2(z)\le 1.\label{49c}
\end{gather}
\label{m23}
\end{lem}

\begin{proof} We just use Lemma~\ref{m9} and Corollary~\ref{m17} to deduce 
\eqref{49b}--\eqref{49d}. Inequality \eqref{49c} is evident.
\end{proof}

Next, we define $U_z$ as in \eqref{47a} using $F_{z,R}$ from
\eqref{sp5}. If $\Gamma=\{z_n\}_{n\ge 1}$, $a_n\in\CC$, $n\ge 1$, and
\begin{equation}
\sum_{n\ge 1} |a_n|^pe^{-ph(z_n)}\rho(z_n)^2\le 1,
\label{h101}
\end{equation}
then for $|z_n|\ge \eta=\max(\eta(R),\eta_1(R))$ 
($\eta(R)$ is introduced in Proposition~\ref{m14})
we define $V_n$ by \eqref{47b}.
Set $\gamma_n=|a_n|e^{-h(z_n)}$.
As above, we obtain
\begin{gather*}
V_n(z_n)=a_n,\\
|V_n(z_k)|e^{-h(z_k)}\le c_1\gamma_n\cdot e^{-cR^2},\qquad |z_k-z_n|\le R\rho(z_n),
\quad k\ne n,\\
\int_{\mathcal D(z,R\rho(z_n))}|V_n(z)|^pe^{-ph(z)}
dm_2(z)\le c_1\gamma^p_n\cdot \rho(z_n)^2,\\
|V_n(z)|e^{-h(z)}\le c_1\gamma_n\cdot
\Bigl(\frac{R^2\min[\rho(z_n),\rho(z)]^2}{e|z-z_n|^2}\Bigr)^{\ve
R^2/5},\,\, |z-z_n|> R\rho(z_n),
\end{gather*}
for some $c,c_1$ independent of $z_n,z,R$.

For $z,\zeta\in\D$ we set
$$
W(z,\zeta)=
\Bigl(\frac{R^2\min[\rho(\zeta),\rho(z)]^2}{e|z-\zeta|^2}\Bigr)^{\ve
R^2/5}(1-\chi_{\mathcal D(\zeta,R\rho(\zeta))}(z)).
$$

Now, to complete the proof as in part (B) of Theorem~\ref{tti01}, we need only to verify that
for every $\delta>0$ there exists arbitrarily large $R$ such that
\begin{align*}
\sum_{k\ge 1}\Bigl(\sum_{n\ge 1}\gamma_n W(z_k,z_n)\Bigr)^p\rho(z_k)^2\le \delta,
\\
\int_{\D}\Bigl(\sum_{n\ge 1}\gamma_n W(z,z_n)\Bigr)^pdm_2(z)\le 1.
\end{align*}
Furthermore, using \eqref{01}, we can deduce these
inequalities from the inequality
\begin{equation}
\int_{\D}\Bigl(\sum_{n\ge 1}\gamma_n W(z,z_n)\Bigr)^pdm_2(z)\le
\delta_1,
\label{f87}
\end{equation}
with small $\delta_1$.
By \eqref{fd144}, for any small $\delta_2>0$ we can find large $R$ such that 
\begin{equation}
\sum_{n\ge 1}W(z,z_n)\le \delta_2,\qquad z\in\D.
\label{h104}
\end{equation}
An estimate analogous to \eqref{fd144} gives us for large $R$
$$
\int_{\D}W(z,z_n)\,dm_2(z)\le c\rho(z_n)^2.
$$
Therefore,
\begin{multline*}
\int_{\D}\Bigl(\sum_{n\ge 1}\gamma_n W(z,z_n)\Bigr)^pdm_2(z)\\
\le \int_{\D}\Bigl(\sum_{n\ge 1}W(z,z_n)\Bigr)^{p-1}\cdot
\Bigl(\sum_{n\ge 1}\gamma^p_n W(z,z_n)\Bigr)\,dm_2(z)\\
\le \delta^{p-1}_2\cdot \sum_{n\ge 1} \gamma^p_n \int_{\D}
W(z,z_n)\,dm_2(z)\le c\delta^{p-1}_2\cdot \sum_{n\ge 1} \gamma^p_n
\rho(z_n)^2 =c\delta^{p-1}_2.
\end{multline*}
This completes the proof of our assertion in the case $p>1$.
If $p=1$, then \eqref{f87} follows from \eqref{h101} and \eqref{h104}.

\end{proof}
\smallskip

\section{The plane case}
\label{spl}

The results of Sections~\ref{sp}--\ref{sin} easily extend to the plane case.
First, we can approximate $h$ by $\log|f|$ for a special infinite product $f$.

\begin{prop} There exist sequences $\{r_k\}$, $\{s_k\}$, $0=r_0<s_0<r_1<\ldots
r_{k}<s_{k}<r_{k+1}<\ldots <\infty$, and a
sequence $N_k$, $k\ge 0$, of natural numbers, such that $N_{k+1}\ge N_{k}$ for large $k$, and
\begin{itemize}
\item[(i)] \quad $\displaystyle \lim_{k\to\infty}\frac{r_{k+1}-r_k}{\rho(r_k)}=\sqrt{2\pi}$,
\quad  $\displaystyle \lim_{k\to\infty}\frac{N_k(r_{k+1}-r_k)}{r_k}=2\pi$,\newline\noindent
\phantom{A} $\displaystyle \lim_{k\to\infty}\frac {r_{k+1}-r_k}{r_{k}-r_{k-1}}=1$,
 \quad $\displaystyle \lim_{k\to\infty}\frac {r_{k+1}-s_k}{r_{k+1}-r_{k}}=\frac12$,
\item[(ii)] \quad if $\Lambda=\bigl\{s_ke^{2\pi i m/N_k}\bigr\}_{k\ge 0,\, 0\le m<N_k}$, and if
\begin{equation*}
f(z)=\lim_{r\to \infty}
\prod_{\lambda\in \Lambda\cap r\D}
\Bigl(1-\frac{z}{\lambda}\Bigr),
\end{equation*}
then the products in the right hand side converge uniformly on compact subsets of the plane, and
\begin{equation*}
|f(z)|\asymp e^{h(z)}\frac{\dist(z,\Lambda)}{\rho(z)},\qquad z\in\CC.
\end{equation*}
\end{itemize}
\label{ll45}
\end{prop}

The proof in the case \condc\ is analogous to that of Proposition~\ref{ll6} in the case
\condd. We need only to mention that
in the estimate \eqref{h4} we use that
\begin{multline*}
\int_0^{y-1} e^{-c_2(r-x)/\rho(x)}\frac{dx}{\rho(x)}\le 
c\int_0^{y-1} e^{-c_2(r-x)/\rho(x)}\frac{dx}{\rho(x)^2}\\
\le c_1\int_0^{y-1}\frac{dx}{(r-x)^2}\le c_1.
\end{multline*}

In the case \conoc, to estimate $\sum U_m$,
$$
U_m=\log|1-s_m^{N_m}z^{-N_m}|,\qquad 0\le m<k,
$$
we divide $m$, $0\le m<k$, into the groups
$$
S_t=\{m:2^t\rho(r)\le \rho(s_m)<2^{t+1}\rho(r)\}, \qquad t\ge 0.
$$

Then for some $M<\infty$
$$
\card S_t\le \frac{cr}{2^{t/M}\cdot 2^t\rho(r)},\qquad r\ge 1,
$$
and
$$ 
\sum_{m\in S_t}e^{-N_m(r-s_m)}\le c_1e^{-cr^2/(2^t\rho(r))}\cdot
\frac{r}{2^{t/M}\cdot 2^t\rho(r)}
\le c_2 2^{-t/M}, \,\, t>0,\, r\ge 1,
$$ 
for some $c,c_1,c_2$ independent of $t>0$, $r\ge 1$,
and
$$
\sum_{m\in S_0}e^{-N_m(r-s_m)}\le \sum_{k\ge 0}e^{-ckr\rho(r)/\rho(r)}\le c_1,
$$
for some $c,c_1$ independent of $r\ge 1$.
\smallskip

Using Proposition~\ref{ll45}, we arrive at analogs of Propositions~\ref{ll7} and \ref{m14}.
For example, we have 

\begin{prop} Given $R\ge 100$, there exists $\eta(R)<\infty$ such that for every $z\in\CC$ with
$|z|\ge \eta(R)$, there exists a function $g=g_{z,R}$ analytic
in $\CC$ such that
\begin{align*}
|g(w)|e^{-h(w)}&\asymp e^{-|z-w|^2/[4\rho(z)^2]},\qquad
w\in \mathcal D(z,R\rho(z)),\\
|g(w)|e^{-h(w)}&\le
c(h)\Bigl[\frac{R^2\min[\rho(z),\rho(w)]^2}{e|z-w|^2}\Bigr]^{R^2/4},\quad
w\in\CC\setminus \mathcal D(z,R\rho(z)).
\end{align*}
\end{prop}
 
After that, the arguments in Sections~\ref{ss}--\ref{sin} extend to the plane case,
and we obtain Theorems~\ref{tt02}--\ref{tti02b}.
\smallskip

\section{Subspaces of large index}
\label{sindex}

Let $X$ be a Banach space of analytic functions in the unit disc, 
such that the operator $M_z$ of multiplication by the independent 
variable $f\mapsto zf$ acts continuously on $X$. Examples
of such spaces are the Hardy spaces $H^p$ and the weighted Bergman spaces $\ahd$, $\ahdp$,
$1\le p<\infty$. 
A closed proper subspace
$E$ of $X$ is said to be $z$-invariant if $M_zE\subset E$.
The index of a $z$-invariant subspace $E$ is defined as
$$
\ind E=\dim E/M_zE.
$$

Every $z$-invariant subspace of the Hardy space $H^2$ has index $1$. 
In 1985, C.~Apostol, H.~Bercovici, C.~Foia{\c s} and C.~Pearcy \cite{BA}
proved (in a non-constructive way) that every space $\ahdd$ has 
$z$-invariant subspaces of index equal to $1,2,\ldots,+\infty$.
Later on, H.~Hedenmalm \cite{HE} and H.~Hedenmalm, S.~Richter, K.~Seip \cite{HRS}
produced concrete examples of $z$-invariant subspaces of index bigger than $1$
in $\mathcal{A}^2_{0}(\D)$, $\mathcal{A}^p_{h_\alpha}(\D)$, with 
$h_\alpha(z)=\alpha\log \frac{1}{1-|z|}$, $\alpha\ge 0$. These examples are based
on Seip's description of sampling sets in $\mathcal{A}^p_{h_\alpha}(\D)$.

In this section we give a construction of $z$-invariant subspaces of large index in $\ahdp$,
$1\le p<\infty$, 
with $h$ satisfying \conod\ or \condd\ based on our results above. For other constructions
suitable for large classes of Banach spaces of analytic functions in the unit disc
and for other information on index of $z$-invariant subspaces see \cite{BOR}, \cite{AB},
\cite{BHV}.

Given a non-empty subset $\Lambda\subset\D$, $1\le p<\infty$, set
$$
I(\Lambda)=\{f\in\ahdp:f(\lambda)=0,\,\lambda\in\Lambda\}.
$$
If $I(\Lambda)\ne \{0\}$, then $I(\Lambda)$ is a closed $z$-invariant subspace
of $\ahdp$, and $\ind I(\Lambda)=1$.

Given $z$-invariant subspaces $E_\alpha$, $\alpha\in A$, denote by 
$\vee_{\alpha\in A}E_\alpha$ the minimal $z$-invariant subspace containing all $E_\alpha$.
It is known that if $\ind E_\alpha=1$, $\alpha\in A$, then
$$
\ind \vee_{\alpha\in A}E_\alpha\le \card A.
$$

\begin{thm} If $h$ satisfies either {\rm\conod} or {\rm\condd}, and if $1\le p<\infty$, then
there exist subsets $\Lambda_d\subset\D$, $0\le d<\infty$, such that
$$
\ind \vee_{0\le d< u}I(\Lambda_d)=u,\qquad 1\le u\le \infty.
$$
\end{thm}

\begin{proof}
We restrict ourselves by the (most difficult) case $u=+\infty$, and use the method proposed in \cite{HRS}.
First, by \eqref{j41}, we can find $\tilde h>h$ satisfying the same conditions as $h$ and such that
\begin{gather*}
\tilde h(r)=(1+o(1))h(r), \qquad r \to 1,\\
\tilde \rho(r)=(1+o(1))\rho(r), \qquad r \to 1,\\
\log \frac 1{\rho(r)}=o(\tilde h(r)-h(r)), \qquad r \to 1,
\end{gather*}
where $\tilde \rho(r)=\bigl[(\Delta \tilde h)(r)\bigr]^{-1/2}$, $0\le r<1$.
Then we apply Proposition~\ref{ll6} to $\tilde h$ to obtain 
$\Lambda=\bigl\{s_ke^{2\pi i m/N_k}\bigr\}_{k\ge 0,\, 0\le m<N_k}$, and 
$f\in \mathcal{A}_{\tilde h}(\D)$ such that $f(0)=1$,
\begin{gather*}
|f(z)|\asymp e^{\tilde h(z)}\frac{\dist(z,\Lambda)}{\tilde \rho(z)},\qquad z\in\D,\\
|f'(\lambda)|\asymp 
\frac{e^{\tilde h(\lambda)}}{\tilde \rho(\lambda)},\qquad \lambda\in\Lambda.
\end{gather*}
An argument similar to that in the proof of Lemma~\ref{m9} shows that for every $g\in\ahdp$
we have
\begin{equation}
\lim_{|z|\to 1}|g(z)|e^{-\tilde h(z)}=0.
\label{h106}
\end{equation}
For $g\in\ahdp$ and for $k\ge 1$ by the residue calculus we have
$$
|g(0)|=\frac{|g(0)|}{|f(0)|}\le \sum_{\lambda\in\Lambda\cap \mathcal D(r_k)}
\frac{|g(\lambda)|}{|\lambda f'(\lambda)|}+
\frac{1}{2\pi}\int_{r_k\mathbb T}\frac{|g(\zeta)|}{|\zeta f(\zeta)|}\,|d\zeta|.
$$
Passing to the limit $k\to\infty$, and using \eqref{h106} and the fact that $\Lambda$ is $d_\rho$-separated
and hence, 
$$
\sum_{\lambda\in\Lambda} \tilde\rho(\lambda)^2 \le c,
$$
we conclude that
\begin{multline}
|g(0)|\le \sum_{\lambda\in\Lambda}\frac{|g(\lambda)|}{|\lambda f'(\lambda)|}\\ 
\le c \Bigl(\sum_{\lambda\in\Lambda}|g(\lambda)|^p 
e^{-p\tilde h(\lambda)}\tilde\rho(\lambda)^2 \Bigr)^{1/p}\Bigl(\sum_{\lambda\in\Lambda}  
\tilde\rho(\lambda)^2 \Bigr)^{(p-1)/p}\\  \le
c\Bigl(\sum_{\lambda\in\Lambda}|g(\lambda)|^p e^{-ph(\lambda)}\rho(\lambda)^2\Bigr)^{1/p}=
c\|g\|_{p,h,\Lambda}.\label{k48}
\end{multline}

We should note here that $\Lambda$ is neither a sampling set
for $\ahdp$ nor that for $\mathcal{A}_{\tilde h}(\D)$.

Now we set 
\begin{gather*}
\Lambda^*_d=\bigl\{s_ke^{2\pi i m/N_k},\,k=2^{d+1}(2v+1),\, v\ge 0,\,0\le m<N_k \bigr\},\\
\Lambda_d=\Lambda\setminus\Lambda^*_d,\qquad d\ge 0.
\end{gather*}
It is clear that
\begin{gather*}
D^{+}_\rho(\Lambda,\D)=D^{-}_\rho(\Lambda,\D)=\frac 12,\\
D^{+}_\rho(\Lambda_d,\D)=D^{-}_\rho(\Lambda_d,\D)=\frac {1-2^{-d-1}}2,\qquad d\ge 0.
\end{gather*}
Using Theorem~\ref{tti01}, we obtain that $I(\Lambda_d)\ne\{0\}$, $d\ge 0$.
To complete our proof we apply the following criterion from \cite[Theorem 2.1]{HRS}.

{\it Suppose that for every $d\ge 0$ there exists $c_d>0$ such that
\begin{equation}
c_d|g(0)|\le \|g+g_1\|_{p,h},\qquad g\in I(\Lambda_d),\,g_1\in 
\vee_{l\ge 0,\, l\ne d}I(\Lambda_l).
\label{k49}
\end{equation}
Then
$$
\ind \vee_{d\ge 0}I(\Lambda_d)=+\infty.
$$
}

It remains to verify \eqref{k49}. Since $\Lambda$ is $d_\rho$-separated
and $\Lambda^*_d$ are pairwise disjoint, by Lemma~\ref{m4} and by \eqref{k48} we obtain
$$
\|g+g_1\|_{p,h}\ge c\|g+g_1\|_{p,h,\Lambda}\ge
c\|g\|_{p,h,\Lambda^*_d}=c\|g\|_{p,h,\Lambda}\ge c_1|g(0)|.
$$
This proves our theorem.
\end{proof}

\bigskip
\bigskip

\noindent \textsc{Alexander Borichev, Department of Mathematics,
University of Bordeaux I, 351, cours de la Lib\'eration, 33405 Talence, France}

\noindent\textsl{E-mail}: \texttt{borichev@math.u-bordeaux.fr}
\medskip

\noindent \textsc{Remi Dhuez, LATP--CMI, Universit\'e de Provence,
39 rue F. Joliot-Curie, 13453 Marseille, France}

\noindent\textsl{E-mail}: \texttt{dhuez@cmi.univ-mrs.fr}
\medskip

\noindent \textsc{Karim Kellay, LATP--CMI, Universit\'e de Provence,
39 rue F. Joliot-Curie, 13453 Marseille, France}

\noindent\textsl{E-mail}: \texttt{kellay@cmi.univ-mrs.fr}

\end{document}